% AMSfnt10.STY (Arabia Septembre 98)
\font\teneurb=eurb10
\font\seveneurb=eurb7
\font\fiveeurb=eurb5
\newfam\eurbfam
\textfont\eurbfam=\teneurb
\scriptfont\eurbfam=\seveneurb
\scriptscriptfont\eurbfam=\fiveeurb

\font\tenmsb=msbm10
\font\sevenmsb=msbm7
\font\fivemsb=msbm5
\newfam\msbfam
\textfont\msbfam=\tenmsb
\scriptfont\msbfam=\sevenmsb
\scriptscriptfont\msbfam=\fivemsb

\font\seveneuf=eufm7
\font\fiveeuf=eufm5

\newfam\euffam
%\textfont\euffam=\teneuf
\scriptfont\euffam=\seveneuf
\scriptscriptfont\euffam=\fiveeuf

%\input diagrammes.tex

%\nofiles
									%\nopagenumbers
%\magnification=1200     
%\font\datefont=cmdunh10

\def\doublemap{\mathrel{\null _{\raise.2ex\hbox{$\textstyle\rightarrow$}}
  ^{\ellower.2ex\hbox{$\textstyle\rightarrow$}}}}

\def\date{le\space\the\day \ifcase\month\or janvier \or f\'evrier\or mars\or
avril\or mai\or juin\or juillet\or ao\^ut\or septembre\or octobre\or
novembre\or d\'ecembre\fi\ {\oldstyle\the\year}}

\def\hfl#1#2{\smash{\mathop{\hbox to 12mm{\rightarrowfill}}
\limits^{\scriptstyle#1}_{\scriptstyle#2}}}

\def\"#1{\if#1i{\accent"7F\i}\else{\accent"7F#1}\fi}
\def\^#1{\if#1i{\accent"5E\i}\else{\accent"5E#1}\fi}

\def\doublemap{\mathrel{\null _{\raise.2ex\hbox{$\textstyle\rightarrow$}}
  ^{\ellower.2ex\hbox{$\textstyle\rightarrow$}}}}

\def\date{le\space\the\day \ifcase\month\or janvier \or f\'evrier\or mars\or
avril\or mai\or juin\or juillet\or ao\^ut\or septembre\or octobre\or
novembre\or d\'ecembre\fi\ {\oldstyle\the\year}}

\def\Z{{\bf Z}}
\def\N{{\bf N}}
\def\C{{\bf C}}
\def\Q{{\bf Q}}

\def\G_m{{\bf G}_m}

\def\e{{\varepsilon}}
\def\X{\chi}

\def\F{{\bf F}}
\let\s\sigma

\let\G\Gamma

\let\w\omega

\def\diag{\mathop{\rm diag}\nolimits}

\def\cite{{\bf [V]}}

\def\lim{\mathop{\rm lim} \nolimits}

\hfuzz=3pt
\overfullrule=0mm
\def\choosenot#1#2{\setbox0=\hbox{$\mathsurround0pt#1#2$}\rlap{\hbox
to\wd0{\hfil$\mathsurround0pt#1/$\hfil}}\hbox{$\mathsurround0pt#1#2$}}

\def\not#1{\mathchoice{\choosenot{\displaystyle}{#1}}{\choosenot{\textstyle}{#1}}{\choosenot{\scriptstyle}{#1}}{\choosenot{\scriptscriptstyle}{#1}}}

% AMSfnt10.STY (Arabia Septembre 98)
\font\teneurb=eurb10
\font\seveneurb=eurb7
\font\fiveeurb=eurb5
\newfam\eurbfam
\textfont\eurbfam=\teneurb
\scriptfont\eurbfam=\seveneurb
\scriptscriptfont\eurbfam=\fiveeurb

\font\tenmsb=msbm10
\font\sevenmsb=msbm7
\font\fivemsb=msbm5
\newfam\msbfam
\textfont\msbfam=\tenmsb
\scriptfont\msbfam=\sevenmsb
\scriptscriptfont\msbfam=\fivemsb

\font\seveneuf=eufm7
\font\fiveeuf=eufm5

\newfam\euffam
%\textfont\euffam=\teneuf
\scriptfont\euffam=\seveneuf
\scriptscriptfont\euffam=\fiveeuf

%\input diagrammes.tex

%\nofiles
									%\nopagenumbers
%\magnification=1200     
%\font\datefont=cmdunh10

\def\doublemap{\mathrel{\null _{\raise.2ex\hbox{$\textstyle\rightarrow$}}
  ^{\ellower.2ex\hbox{$\textstyle\rightarrow$}}}}

\def\date{le\space\the\day \ifcase\month\or janvier \or f\'evrier\or mars\or
avril\or mai\or juin\or juillet\or ao\^ut\or septembre\or octobre\or
novembre\or d\'ecembre\fi\ {\oldstyle\the\year}}

\def\hfl#1#2{\smash{\mathop{\hbox to 12mm{\rightarrowfill}}
\limits^{\scriptstyle#1}_{\scriptstyle#2}}}

\def\"#1{\if#1i{\accent"7F\i}\else{\accent"7F#1}\fi}
\def\^#1{\if#1i{\accent"5E\i}\else{\accent"5E#1}\fi}

\def\doublemap{\mathrel{\null _{\raise.2ex\hbox{$\textstyle\rightarrow$}}
  ^{\ellower.2ex\hbox{$\textstyle\rightarrow$}}}}

\def\date{le\space\the\day \ifcase\month\or janvier \or f\'evrier\or mars\or
avril\or mai\or juin\or juillet\or ao\^ut\or septembre\or octobre\or
novembre\or d\'ecembre\fi\ {\oldstyle\the\year}}

\def\Z{{\bf Z}}
\def\N{{\bf N}}
\def\C{{\bf C}}
\def\Q{{\bf Q}}

\def\G_m{{\bf G}_m}

\def\e{{\varepsilon}}
\def\X{\chi}

\def\F{{\bf F}}
\let\s\sigma

\let\G\Gamma

\let\w\omega

\def\diag{\mathop{\rm diag}\nolimits}

\def\cite{{\bf [V]}}

\def\lim{\mathop{\rm lim} \nolimits}

\hfuzz=3pt
\overfullrule=0mm
\def\choosenot#1#2{\setbox0=\hbox{$\mathsurround0pt#1#2$}\rlap{\hbox
to\wd0{\hfil$\mathsurround0pt#1/$\hfil}}\hbox{$\mathsurround0pt#1#2$}}

\def\not#1{\mathchoice{\choosenot{\displaystyle}{#1}}{\choosenot{\textstyle}{#1}}{\choosenot{\scriptstyle}{#1}}{\choosenot{\scriptscriptstyle}{#1}}}

%\documentclass[11pt]{article}
%\usepackage{graphicx}
%\usepackage{pdfsync}
%\usepackage{amssymb}
%\magnification=1200
%\begin{document} 
\centerline {\bf Alg\`ebres de Hecke affines g\'en\'eriques \ }

\centerline {Marie-France Vign\'eras}

\bigskip {\bf R\'esum\'e \ } Soit 
  $H$ l'alg\`ebre de Hecke   du groupe de Weyl $W$ d'une
 donn\'ee radicielle  bas\'ee $(X,X^\vee,R,R^{\vee},B)$, et d'un poids
g\'en\'erique  $(q_w)_{w\in W}$. 
 Nous   montrerons que $H$ est un  module de type fini
sur son centre, et que
le centre est une
alg\`ebre
\`a engendrement fini. Ceci  
\'etait   connu   apr\`es inversion du poids, mais il
est essentiel de ne pas inverser le poids, dans
l'\'etude
 des modules  des alg\`ebres de Hecke affines  qui  apparaissent
naturellement dans  la th\'eorie
  des repr\'esentations des groupes r\'eductifs $p$-adiques
sur un corps $p$-adique ou de caract\'eristique
$p$. Des applications de ces th\'eor\`emes de finitude \`a la
th\'eorie des modules sont donn\'es dans une seconde partie.  Dans une
troisi\`eme partie, on explicitera les r\'esultats   pour
le groupe
$GL(n)$, et l'on introduira pour ce groupe, les modules supersinguliers.

\bigskip {\bf Introduction \ } 

{\sl Chapitre 1}.   L'alg\`ebre
de Hecke g\'en\'erique $H$ est un module libre sur l'alg\`ebre de
polyn\^omes g\'en\'erique $\Z[q_w, w\in W]$ (not\'ee $\Z[q_*]$ dans la
suite),   ayant une base naturelle
$(T_w)_{w\in W}$,
appel\'ee la ``base de Iwahori-Matsumoto''.  Le groupe de Weyl fini
$W_o$ du syst\`eme de racines $R$   agit naturellement sur $X$. Le
groupe de Weyl $W $ de la donn\'ee radicielle est le produit
semi-direct  
$W=W_o X$.  Les
\'el\'ements
$T_w$  sont inversibles  dans l'alg\`ebre   $H[q_*^{-1}] $ obtenue apr\`es
inversion des poids $q_w, w\in W$. Suivant Bernstein et Lusztig, il existe
une application
$x\to \tilde \theta_x $ de  $X$ dans  l'alg\`ebre de Hecke  
$H[q_*^{-1/2}]$ obtenue apr\`es addition des
inverses des racines carr\'ees des poids,  identifiant l'alg\`ebre de
groupe   $\Z[q_*^{\pm 1/2}][X]$ \`a une sous-alg\`ebre commutative $  A$
de $H[q_*^{-1/2}]$. Pour un \'el\'ement $w=(w_o,x)$ de
$W$, nous posons $E_w:=q_w^{1/2} T_{w_o}\tilde
\theta_x$.   Nous montrons, sans utiliser la th\'eorie
de Bernstein-Lusztig,  en
\'etudiant le d\'eveloppement de
$T_wT_{v^{-1}}^{-1}$ sur la base de Iwahori-Matsumoto pour tout
$w,v\in W$, le r\'esultat fondamental suivant:
 
\bigskip {\bf Th\'eor\`eme 1 \ } {\sl  
$(E_w)_{w\in W}$ est une base   de l'alg\`ebre de Hecke g\'en\'erique $H$
sur $\Z[q_*]$,
  et la matrice de passage avec la base de
Iwahori-Matsumoto $(T_w)_{w\in W}$ est triangulaire, de
coefficients     diagonaux 
$1$, pour l'ordre de Chevalley-Bruhat.}

\bigskip 
La  base 
   $(E_w)_{w\in W}$ de $H$  sur $\Z[q_*]$ a des propri\'et\'es
remarquables qui permettent de montrer  que:

- {\sl  $H\cap  A$   
est un $\Z[q_*]$-module libre de base 
$(E_x=q_x^{1/2}\tilde \theta_x )_{x\in X}$, 

- $H$  est un $H\cap  
A$-module \`a droite de type fini.}

\bigskip   Comme
  $x\mapsto q_x$ est constante sur toute 
$W_o$-orbite de
$X$, l'alg\`ebre commutative $H\cap A$
est stable par l'action naturelle de
$W_o$ sur $A$, et 
  le centre de
$H$ est
$(H\cap   A )^{W_o}$. Nous obtiendrons  le  th\'eor\`eme
structural suivant:
 
\bigskip {\bf Th\'eor\`eme 2 \ } {\sl L'alg\`ebre de Hecke g\'en\'erique
$H$ est un module de type fini sur son centre $Z(H)=(H\cap   A
)^{W_o}$,   le centre et $H\cap A$ sont des alg\`ebres commutatives de
type fini, et 
$$H=\oplus_{w_o\in W_o}T_{w_o} J(w_o)$$
pour des id\'eaux fractionnaires $J(w_o)$ de $A\cap H$.}

\bigskip En g\'en\'eral les   $A\cap H$-modules $J(w_o)$ ne 
sont plats,   i.e.  projectifs puisque $A\cap
H$ est un anneau noetherien. Des contre-exemples ont \'et\'e donn\'es par Rachel Ollivier 
pour  les alg\`ebres associ\'ees au groupe $GL(3)$.

\bigskip  {\sl Chapitre 2}. Il contient  des applications 
des th\'eor\`emes de finitude du chapitre 1 \`a la th\'eorie des 
modules. On fixe un morphisme d'anneau commutatif  $\phi:\Z[q_*]\to R$
et l'on note
$H_{\phi}$ l'alg\`ebre de Hecke sp\'ecialis\'ee. Pour tout
morphisme d'anneau
$\X:A\cap H
\to R$   prolongeant
$\phi$, le module induit de $\X$ 
$$I (\X):= H \otimes _{ A\cap H , \X} R$$  
(le $H$-module universel engendr\'e par le caract\`ere $\X$ de $A\cap
H$) est appel\'e   un module standard (\`a gauche) de $H_{\phi}$. Si 
les images des poids $q_s, \ s\in S, $
 sont inversibles et ont des racines carr\'ees dans
$R$ (on peut se passer des racines carr\'ees), c'est simplement le
module   induit de l'unique  morphisme d'anneau $\X_A:A\to R$
prolongeant $\X$,
$$I(\X_A)=H[q_*^{-1/2}]\otimes_{A, \X_A} R.$$
 Le th\'eor\`eme 2 implique que:

-    $I (\X)$
  est un $R$-module  de type fini. 

- Lorsque  
$R$ est un corps alg\'ebriquement clos,  un
$H_{
\phi}$-module simple est de dimension finie si et seulement si  
$Z(H)$ agit  par un caract\`ere, appel\'e le caract\`ere
central.

\bigskip  {\bf Th\'eor\`eme 3 \ } {\sl Si
$R$ est un corps alg\'ebriquement clos, un $H_{\phi}$-module simple
de dimension finie
est   quotient  d'un  $H_{\phi}$-module  standard.}

\bigskip Dans le cas classique, le th\'eor\`eme est aussi vrai avec
{\sl sous-module} au lieu de  quotient. Je ne sais pas s'il reste vrai
dans  le cas  o\`u $\phi(q_s)=0$   pour tout $s\in S$ de longueur non
nulle, apparaissant dans l'\'etude des repr\'esentations modulo $p$ des groupes r\'eductifs
$p$-adiques. 

\bigskip
Soit $\overline \Q_p$    une cl\^oture alg\'ebrique du corps des
nombres $p$-adiques, $\overline \Z_p$  
 son anneau d'entiers  de corps r\'esiduel $\overline
\F_p$. 
On note $i_p:\overline \Z_p\to \overline \Q_p$ l'inclusion et
$r_p:\overline \Z_p\to \overline \F_p$ la r\'eduction. On fixe
un morphisme d'anneau $\phi:
\Z[q_*]\to
\overline \Z_p$  tel que $\phi(q_s)\neq 0$   pour tout $s\in S$.  La
finitude de l'alg\`ebre de Hecke g\'en\'erique  $H$ 
permettra de transf\'erer certaines propri\'et\'es de la
th\'eorie classique des $H_{i_p\phi}$-modules \`a celle  des
$H_{r_p\phi}$-modules, par r\'eduction.

\bigskip Un $H_{i_p\phi}$-module
$V$ sera appel\'e entier s'il est engendr\'e sur $\overline \Q_{p}$ par un sous-$\overline \Z_p$-module libre $L$ stable par 
$H_{\phi}$, appel\'e une structure
enti\`ere de $V$.
Pour un 
$H_{i_p\phi}$-module simple de dimension finie, la propri\'et\'e d'\^etre
entier se voit sur le caract\`ere central.

\bigskip {\bf Th\'eor\`eme 4   \ } {\sl Un  $H_{i_p\phi}$-module  
 simple  de dimension finie est entier si et seulement si son
caract\`ere central est entier.} 

\bigskip  Tout caract\`ere  $\X:A\cap H\to \overline \Q_p$   prolongeant $\phi$ s'\'etend  en un caract\`ere $\X_A$ de $A$ et  s'identifie \`a un  caract\`ere $ \X_{X}$ de $X$   en posant $\X_A(\tilde \theta_x)= \X_{X} (x)$ pour tout $x\in X$.  L'isomorphisme   de Bernstein-Lusztig
$\tilde \theta: \Z[q_*^{\pm 1/2}][X]\simeq A$ n'est pas compatible aux   structures enti\`eres naturelles puisque une $\Z[q_*]$-base de $A\cap H$ est $(E_x= q_x^{1/2}\tilde \theta_x)_{x\in X}$. 

\bigskip {\bf Th\'eor\`eme 4 (suite) \ } {\sl Soit $ \X_{X}:X \to \overline \Q_p$   un caract\`ere de $X$.  Le
$H_{i_p\phi}$-module  standard $I( \X_A) $  associ\'e  est entier si et seulement si $ \X_{X} (x)
\phi(q_x)^{-1/2}\in \overline \Z_p $  pour tout $x\in X$.}
 
\bigskip Comme $q_x$ est constant sur une $W_o$-orbite de $X$, la
condition  est $W_o$-invariante.  Le module standard est  de dimension $|W_{o}|$ sur  $\overline \Q_p$ et admet une base canonique:   l'image canonique de $(T_{w_{o}})_{w_{o}\in W_{o}}$.

\bigskip {\bf Th\'eor\`eme 5   \ } {\sl   
Soit $\X :A\cap H\to \overline \Z_p$  un morphisme d'anneau prolongeant $\phi$. Alors

(i)  Le module standard
 $I(\X,\overline \Q_{p}):= H \otimes_{i_{p}\X}\overline \Q_{p}=I(\X) \otimes_{i_{p}}\overline \Q_{p}$  admet une structure enti\`ere canonique  contenant la base canonique, isomorphe au quotient de $I(\X)$ par son sous-groupe de torsion $I(\X)_{tor}$.

 (ii)  Les propri\'et\'es  suivantes sont  \'equivalentes:
 
a) Le module standard $I(\X,\overline \F_{p}):= H\otimes_{r _{p}\X}\overline \F_{p}=I(\X)\otimes_{r_{p}}\overline \F_{p}$ est de dimension $|W_{o}|$ sur $\overline \F_{p}$,
 
 b) La r\'eduction de la structure enti\`ere canonique $M(\X):=I(\X)/I(\X)_{tor}$ de $I(\X,\overline \Q_{p})$ est isomorphe au 
 module standard $I(\X,\overline \F_{p})$,
 
 c) Le $\overline \Z_{p}$-module $I(\X)$ est sans torsion.

Elles sont v\'erifi\'ees si le $A\cap H$-module \`a droite $H$ est plat. }
\bigskip Lorsque $r_{p}\phi(q_{s})= 0$ pour $s\in S$,  les propri\'et\'es (ii) sont   v\'erifi\'ees pour $GL(2)$ mais pas pour $GL(3)$, d'apr\`es les r\'esultats de Rachel Ollivier.  Donc en g\'en\'eral, le
$A\cap H$-module
\`a droite   de type fini $H$ n'est pas
plat.   Il serait int\'eressant de savoir s'il existe un caract\`ere $\X':A\cap H \to \overline \Z_{p}$ de m\^eme r\'eduction que $\X$   mais tel que la  r\'eduction  $r_{p}M(\X')$ n'ait pas la m\^eme semi-simplifi\'ee que $r_{p}M(\X)$,  et s'il existe un caract\`ere $\X$ tel que le noyau du morphisme canonique   $$I(\X,\overline \F_{p})\to \prod_{\X', r_{p}(\X')=r_{p}(X)}r_{p}M(\X')$$ n'est pas nul.

\bigskip {\sl Chapitre 3}. Il contient une \'etude  d\'etaill\'ee de
l'alg\`ebre de Hecke g\'en\'erique du groupe  
lin\'eaire $GL(n)$, et  l'introduction de ses modules
supersinguliers. 

Le groupe 
$X\simeq \Z^n$ s'identifie au groupe $Mor (GL(1), GL(n))$  des
cocaract\`eres sur lequel le groupe de Weyl fini $W_o$ identifi\'e \`a
$S_n$, agit naturellement. Le groupe de Weyl $W$ est isomorphe au
produit semi-direct naturel 
$S_n
\Z^n$.  On a 
$q _s=q$ pour tout $s\in S$. 
L'alg\`ebre de Hecke
g\'en\'erique $H$ est une   $\Z[q]$-alg\`ebre. 
On donnera (3.1-2) une description explicite de $A\cap H$ et de ses
caract\`eres dont on  en d\'eduira les r\'esultats suivants. On fixe
une sp\'ecialisation $\phi:\Z[q]\to \overline \Z_p$ telle que
$\phi(q)\neq 0$
  n'est pas une unit\'e $p$-adique.

\bigskip  {\bf Th\'eor\`eme 6 \ }  
{\sl  Pour $GL(n)$,
 tout caract\`ere
$\X:A\cap H
\to
\overline \F_p$  nul sur $q$, se rel\`eve en  un caract\`ere
$\tilde
\X:A\cap H
\to
\overline
\Z_p$  tel que $\tilde \X(q)=\phi(q)$.}

\bigskip 
Il existe des $\overline \F_p$-repr\'esentations irr\'eductibles de
$GL(2,\Q_p)$ qui ne sont pas des sous-quotients d'induites paraboliques
propres, appel\'ees supersinguli\`eres.
Par le foncteur des invariants par le sous-groupe d'Iwahori, elles sont
en bijection avec des $H_{r_p\phi}$-modules simples standards. Une
g\'en\'eralisation naturelle de ces modules  est celle-ci:

\bigskip {\bf D\'efinition   \ } {\sl Pour
$GL(n)$, on
dira qu'un 
$H_{r_p\phi}$-module ayant un caract\`ere central $\w$ est 
supersingulier, si  $\w$ se prolonge en un (unique) caract\`ere de $A\cap
H$ fixe par $S_n$.}

\bigskip Pour un caract\`ere de
$  A\cap H $  nul sur
$q$, {\sl \^etre  fixe par $S_n$ signifie \^etre   nul   sur
$E_x$ pour tout  
$x\in X$ qui n'est pas une puissance du plongement diagonal $\delta$ de
$GL(1)$ dans $GL(n)$}. La valeur en $E_{\delta}$ est toujours inversible car
$E_\delta=T_\delta$ est inversible dans $H$. Pour
$z\in
\overline \F_p^*$,  le caract\`ere 
$  A\cap H\to \overline \F_p$  nul sur
$q$, fixe par $S_n$, tel que $\X(E_{\delta})=z$
sera not\'e  
$\X_z$, et sa restriction  au centre sera not\'ee
$\w_z$. Le $H_{r_p\phi}$-module standard $I(\X_z)$ est l'unique module
standard supersingulier de caract\`ere central $\w_z$.

\bigskip {\bf  Proposition 7 \ } {\sl Pour
$GL(n)$, un module  standard supersingulier $I(\X_z)$  
 est de longueur $\geq 2^{n-2}$. Tout $H_{r_p\phi}$-module simple
supersingulier de caract\`ere central $\w_z$ est quotient de $I(\X_z)$.
Pour
$n=2 $, $I(\X_z)$ est irr\'eductible de dimension $2$. }

\bigskip Le cas $n=2$ est trait\'e dans un article pr\'ec\'edent [Vign\'eras].  

\bigskip Les outils actuels de la th\'eorie des repr\'esentations
s'adaptent mal aux repr\'esentations lisses modulo
$p$ des groupes r\'eductifs $p$-adiques pour lesquelles on ne sait pour
ainsi dire rien. Ce travail est le premier pas vers une
classification des modules simples    modulo
$p$ de leurs alg\`ebres de Hecke-Iwahori, probl\`eme qui semble plus
accessible et qui est reli\'e au pr\'ec\'edent par le foncteur des
invariants par un sous-groupe d'Iwahori. La situation est tr\`es
diff\'erente de celle des groupes r\'eductifs  finis sur un corps de
caract\'eristique
$p$, o\`u tous ces probl\`emes sont r\'esolus: les modules simples
modulo
$p$ sont class\'es   par la th\'eorie du plus haut poids, le
foncteur des invariants par un sous-groupe de Borel respecte
l'irr\'eductibilit\'e l\`a o\`u il ne s'annule pas, et les modules
simples modulo
$p$ des
alg\`ebres de Hecke-Iwahori  sont tous de dimension
$1$.

\bigskip L'auteur remercie Rachel Ollivier dont les
travaux sur les modules standards de type $A_2$  ont confirm\'e la
possibilit\'e d'une th\'eorie int\'eressante des alg\`ebres de Hecke
affines de param\`etre
$q=0$ et qui ont permis de corriger  une version pr\'ec\'edente du th\'eor\`eme 5, Jean-Francois Dat  pour son
\'etude parall\`ele des repr\'esentations
enti\`eres $p$-adiques des groupes r\'eductifs
$p$-adiques, Peter Schneider et les participants de la conf\'erence 
``Fonctorialit\'e de Langlands: progr\`es recents'' (Luminy, CIRM,
21/29 juin 2002)  dans laquelle  ces r\'esultats furent expos\'es,  
l'Institut de Math\'ematiques de Jussieu pour un environnement
de recherche remarquable, et le C.N.R.S. pour son soutien financier.

\vfill\eject \centerline {{\bf Chapitre 1 \ }}

\bigskip Soit  $(X,X^\vee,R,R^{\vee},B)$ une donn\'ee radicielle
bas\'ee r\'eduite [Lusztig1 1 page 600-601].
Le groupe de Weyl  fini $W_o$ est le  syst\`eme de Coxeter
avec $S_o=\{ s_{\alpha} \ | \  \alpha \in B\}$ comme ensemble de
r\'eflexions simples; le produit semi-direct 
$W=W_o.X$ est le groupe de Weyl de la donn\'ee radicielle;  
 le groupe de Weyl affine est $W_{aff}:=W_o.Q(R)$ o\`u $Q(R)$
est le sous-$\Z$-module  de $X$ engendr\'e par $R$. Il existe un
sous-groupe commutatif $\Omega$ de $W$ form\'e 
 d'\'el\'ements tel
que $W_{aff}\cap
\Omega=\{1\}$ et
$W=\Omega W_{aff}$.

Traditionellement la notation pour le produit est additive
dans le
$\Z$-module libre
$X$ et multiplicative dans les groupes
$W_o$ et dans $W$. Pour cette raison on notera
$e^x$ l'\'el\'ement
$x\in X$ plong\'e dans $W$ qui s'\'ecrira dor\'enavant $$W=W_oe^X=
 \Omega W_{aff}.$$  Un
\'el\'ement g\'en\'eral de $W$ s'\'ecrira donc $w=w_oe^x=uw_{aff} $
avec
$w_o\in W_o,
\ x\in X,
\ u\in U, \  w_{aff}\in W_{aff},$ uniques. La valeur en $ e^x$ d'une  
fonction
$w\to f_w$ sur $W$ sera aussi not\'ee $  f_x$.

Soit $R_m$ l'ensemble des racines $\alpha\in R$ telles que
$\alpha^\vee$ est minimal pour l'ordre partiel $\preceq $ d\'efini sur
les coracines $R^\vee$  par: $\alpha_1^\vee \preceq \alpha_2^\vee $ si
$\alpha_2^\vee-\alpha_1^\vee \in \sum _{\alpha\in B} \N \alpha^\vee$.
Le groupe de Weyl affine
$W_{aff}$ est un syst\`eme de Coxeter avec $$S=S_o \cup \{
s_{\alpha}e^\alpha \ | \ 
\alpha \in R_m\}$$ comme ensemble de r\'eflexions simples.  Le groupe
$\Omega$ normalise $S$. On note
$\leq $ l'ordre de Chevalley-Bruhat, $\ell$ la longueur de $W$,  
d\'efinis par extension
   de l'ordre et de la longueur de
$(W_{aff},S)$, et $X_{dom}$ 
le monoide de type fini form\'e par les
\'el\'ements dominants $x$ de
$X$, i.e.   $(x,\alpha^\vee)\geq 0$ pour toute racine
positive $\alpha $. 
Les d\'efinitions et rappels ci-dessus sont d\'etaill\'es
dans
 l'appendice.

\bigskip

 \bigskip \centerline {{\bf Alg\`ebre de Hecke g\'en\'erique \ }}
 \bigskip L'alg\`ebre  de Hecke g\'en\'erique du groupe
de Weyl de la donn\'ee radicielle bas\'ee sera d\'efinie ``\`a la
Iwahori-Matsumoto'' en utilisant la d\'ecomposition en produit
semi-direct
$W=\Omega W_{aff}$. Lorsque la donn\'ee radicielle est torale,
$W=\Omega=X,$   l'alg\`ebre de Hecke g\'en\'erique est simplement
l'alg\`ebre du groupe $\Z[X]$.  

\bigskip {\bf Poids g\'en\'erique \ }  Lorsque la donn\'ee radicielle
bas\'ee est non torale,   soient $(q_s)_{ s\in S} $ des
{\sl ind\'etermin\'ees
 telles que $q_s=q_{s'}$ si $s,s'\in S$ sont conjugu\'es
dans $W$}. 

  La $\Z$-alg\`ebre de
polyn\^omes  engendr\'ee  par ces ind\'etermin\'ees  est not\'ee 
$\Z[q_*]
$. On connait un crit\`ere pour que deux \'el\'ements de $S$ soient
conjugu\'es dans
$W_{aff}$ [Bourbaki GAL IV.1 proposition 3, page 12]. Deux
\'el\'ements de
$S$ peuvent
\^etre conjugu\'es dans
$W$ sans l'\^etre dans $W_{aff}$, et l'alg\`ebre de polyn\^omes
associ\'ee
\`a $W$ est une sp\'ecialisation de celle associ\'ee \`a $ W_{aff} $.

\bigskip Pour tout $n \in \N^S$ on dit que
$q_n:=\prod_{s\in S} q_s^{n(s)}$   
  est un mon\^ome en $q_*$. Si $n$ est la fonction
``multiplicit\'e de
$s$ dans une d\'ecomposition r\'eduite'' de  $w_{aff}\in W_{aff}$, 
alors le fait que
$q_s=q_{s'}$ si
$s,s'$ sont conjugu\'ees dans $W$, montre que le mon\^ome $q_n$ ne
d\'epend pas de la    d\'ecomposition r\'eduite, et sera not\'e
$q_{w_{aff}}$. On ne risque pas de contresens en posant
$q_w=q_{w_{aff}}$ pour tout $w\in
 \Omega w_{aff}
\Omega
$. On a
$q_w=q_{w^{-1}}$. 
On a $q_w=q_{w'}$ pour deux \'el\'ements conjugu\'es
$w,w'$ de
$W$. On notera $q_x=q_{e^x}$ pour $x\in X$. La fonction $x\to q_x$ est
constante sur les $W_o$-orbites de $X$ puisque
$e^{w_o(x)}=w_oe^xw_o^{-1}$.
 
 On dit que $(q_w)_{w\in W}$ est un poids g\'en\'erique
de $W$.

\bigskip {\bf 1.1  \ D\'efinition \ } {\sl  L'alg\`ebre
  de Hecke g\'en\'erique $H$ est la   
 $\Z[q_*]$-alg\`ebre, libre de base $(T_w)_{w\in W}$ comme 
$\Z[q_*]$-module, de produit v\'erifiant les
relations  
 
- $(T_s+1)(T_s-q_s)=0 \  \ \ {\rm pour \ tout \ } \ \ s\in S, $ 

- $T_{ww'}=T_{w}T_{w'}$ si $\ell(ww')=\ell(w)+\ell(w')  \ \ {\rm pour \
tout \ } \ \ w,w'\in W$.}

\bigskip {\bf Remarque \ } Dans [Lusztig1 3.1 page 606], l'alg\`ebre
$\Z[q_*]$   est remplac\'ee par l'anneau des polyn\^omes de Laurent
$ \Z[v ^{\pm 1}]$, que l'on peut voir comme une localisation d'une
sp\'ecialisation
$\phi_L:\Z[q_*]\to \Z[v]$ d\'efinie par
$\phi_L(q_s)=v^{2n_s}$ pour des entiers $(n_s \geq 0)_{s\in S}$ tels
que $n_s=n_{s'}$ si
$s,s'$ sont conjugu\'ees dans $W$. La
$\Z[v]$-alg\`ebre sp\'ecialis\'ee
$H_{\phi _L}$ est moins g\'en\'erale que $H$ mais est suffisante
pour les applications \`a la th\'eorie des repr\'esentations des
groupes r\'eductifs $p$-adiques. La
$\Z[v^{\pm 1}]$-alg\`ebre
$H_{\phi_L}[v^{-1}]$ est consid\'er\'ee dans 
[Lusztig1 3.2 page 607] ou [Lusztig3 3, si $W=W_{aff}$].

\bigskip   Le groupe $\Omega$ se plonge dans $H$ par $u\to T_u$, le
$\Z[q_*]$-module libre de base
$(T_w)_{w\in
W_{aff}}$ est une sous-alg\`ebre $H_{aff}$ de $H$ normalis\'ee par
$T_u, u\in \Omega$, et   $$H\simeq \Z[\Omega]. H_{aff}.$$ L'alg\`ebre
$H_{aff}$ est une sp\'ecialisation de  l'alg\`ebre de  Hecke
g\'en\'erique de
$W_{aff}$.

  Les
\'el\'ements de la base  
$(T_w)_{w\in W}$ deviennent inversibles si l'on inverse les
param\`etres, i.e. $T_w$ est inversible dans la $\Z[q_*^{\pm 1
}]$-alg\`ebre  
$H[q_*^{-1 }] $ qui contient $H$.

\bigskip  {\bf 1.2  \ Lemme fondamental \ } {\sl  Soient $ w,v\in
W$. Alors 

1) $q_{wv}q_w^{-1}q_v=c_{w,v}^2$ est le carr\'e d'un mon\^ome
$c_{w,v}$ en
$q_*$ divisant   $q_v$, 

2) Soient $ u\in
\Omega, \ s_1, \ldots, s_n\in S$ tels  que $v\in W_{aff}u$  et
$vu^{-1}=s_1\ldots s_n $ soit une d\'ecomposition r\'eduite.
Alors
$$c_{w,v}  T_w   T_{v^{-1}}^{-1}
\ =  \ T_{wv}\ + \ \sum 
  \   T_{ws_1\ldots \hat s_{i_1} \ldots 
\hat s_{i_2}
\ldots
\hat s_{i_q}\ldots s_{n}u  }  \ 
  \lambda_{s_{i_1}}  \lambda_{s_{i_2}} \ldots  
\lambda_{s_{i_q}}   \ \prod _{j\in \{1, \ldots, n\}- \{i_1, i_2,
\ldots, i_q\}}
 q_{s_j}^{k_j}$$
o\`u  
$\lambda_s:=1-q_s$ pour $s\in S$, les
$k_j$ sont des entiers $\geq 0$, et la somme est prise 
sur les suites $(s_{i_1},s_{i_2},\ldots, s_{i_q})$ de longueur $\geq
1$   extraites de la suite $(s_1,\ldots, s_n)$ qui sont
$w$-distingu\'ees   au sens
suivant:  

  $\s_{i_1-1}:=ws_1\ldots s_{i_1-1}< \s_{i_1-1}s_{i_1}$

$\s_{i_2-1}:=\s_{i_1-1} s_{i_1+1}\ldots s_{i_2-1}<
\s_{i_2-1}s_{i_2},$

$\ldots$

$\s_{i_q-1}:=\s_{i_{q-1}-1} s_{i_{q-1}+1}\ldots s_{i_q-1}<
\s_{i_q-1}s_{i_q}.$ 

3) $ws_1\ldots \hat s_{i_1} \ldots 
\hat s_{i_2}
\ldots
\hat s_{i_q}\ldots s_{n}u  < wv$ pour toute suite
$(s_{i_1},s_{i_2},\ldots, s_{i_q})$ de longueur $\geq 1$   extraite 
de la suite $(s_1,\ldots, s_n)$, qui est 
$w$-distingu\'ee.}

\bigskip  Il est pratique  
d'introduire de nouvelles ind\'etermin\'ees $(q_s^{1/2})_{ s\in S} $
\'egales si
$s,s'$ sont conjugu\'es dans $W$, dont les carr\'es sont $(q_s)_{ s\in
S}$.  
  Pour tout $s\in S$   on note  
$  \tilde \lambda_s:=q_s^{-1/2}\lambda_s$, et 
pour tout $w\in W$ on note $ \tilde T_w:=q_w^{-1/2}T_w$. Le lemme
fondamental r\'esulte essentiellement des relations 
$$\tilde
 T_s^{-1}=\tilde T_s+\tilde \lambda_s,$$
 $$\tilde T_w \tilde T_s^{-1}=\tilde T_{ws} \ \
{\rm  si
\ }
\ ws<w ,
\ {\rm  et \ } \ \tilde T_w \tilde T_s^{-1}=\tilde T_{ws}+\tilde
\lambda_s\tilde T_w 
\ {\rm  si \ } \ w<ws.
 $$

\bigskip Preuve.  La d\'emonstration de 3)  
est 
   faite
  dans  ([Ha] prop. 5.5).  On y d\'emontre l'in\'egalit\'e au
sens large $ws_1\ldots
\hat s_{i_1}
\ldots 
\hat s_{i_2}
\ldots
\hat s_{i_q}\ldots s_{n}  \leq wv$, l'\'egalit\'e impliquerait en
simplifiant $\hat s_{i_1}
\ldots 
\hat s_{i_2}
\ldots
\hat s_{i_q}= s_{i_1}
\ldots 
 s_{i_q}$ ce qui est faux si $q\neq 0$.  On a donc l'in\'egalit\'e
stricte.

Montrons 1). 

 a) Si $v=s\in S$  on a $q_{ws}=q_w q_s^{\e}$ o\`u
$\ell(ws)=\ell(w)+\e$; donc 1) est v\'erifi\'e pour $(w,s)$ avec 
$c_{w,s}=  q_s$ si $ws>w$ et $ 1$ si $ws<w$.

b) Si $v\in W_{aff}$ et $v=s_1\ldots s_n$ est une
d\'ecomposition r\'eduite,  on a par induction $q_ {wv}=q_w
q_{s_1}^{\e_1}\ldots  q_{s_v}^{\e_n}$ o\`u 
 $\ell(ws_1)=\ell(w)+\e_1, \ \ell(ws_1s_2)=\ell(ws_1)+\e_2,\ \ldots, \
\ell(wv)=\ell(ws_1\ldots s_{n-1})+\e_n.$ On a $q_v=q_{s_1} \ldots 
q_{s_n} $ et 1) est v\'erifi\'e pour $(w,v)$ avec 
$$c_{w,v}= \prod _{1\leq i \leq n, \e_i=1} q_{s_i} $$

c) Si $v=v_{aff}u, v_{aff}\in W_{aff}, u\in \Omega,$ on a
$q_{wv}=q_{wv_{aff}}, q_v=q_{v_{aff}}$. Donc par b), 
le lemme est
v\'erifi\'e pour $(w,v)$ avec $c_{w,v}=c_{w,v_{aff}}$.

 \bigskip Montrons 2). 
Il est imm\'ediat que le lemme fondamental est v\'erifi\'e si $v=s\in
S$.
Soient $s_1, \ldots, s_n\in S$ tels que  $v=s_1\ldots s_n$ soit une
d\'ecomposition r\'eduite. En appliquant la relation 
pr\'ec\'edente plusieurs fois, on obtient
$$ 
\tilde T_w 
\tilde T_{v^{-1}}^{-1}=\tilde T_{wv}+\sum 
  \  \tilde T_{ws_1\ldots \hat s_{i_1} \ldots 
\hat s_{i_2}
\ldots
\hat s_{i_q}\ldots s_{n}  }  \ 
 \tilde  \lambda_{s_{i_1}} \tilde  \lambda_{s_{i_2}} \ldots  \tilde 
\lambda_{s_{i_q}} $$
o\`u  la somme est prise 
sur les suites $(s_{i_1},s_{i_2},\ldots, s_{i_q})$  de longueur $\geq
1$  extraites de la suite $(s_1,\ldots, s_n)$ qui sont
$w$-distingu\'ees. Le lemme fondamental   pour $v\in W_{aff}$ est
alors 
\'equivalent \`a l'assertion suivante
sur les fonctions $q_w$.

\bigskip    {\sl  Soit $w\in W, s_1,\ldots, s_n \in S$ tels que
$s_1\ldots s_n$ soit une d\'ecomposition r\'eduite. Pour une suite
$(s_{i_1},s_{i_2},\ldots, s_{i_q})$   extraite 
de la suite $(s_1,\ldots, s_n)$ et $w$-distingu\'ee, on a 

$$q_{ws_1\ldots s_n}\ =  \ 
q_{ws_1\ldots \hat s_{i_1} \ldots \hat s_{i_2} \ldots
\hat s_{i_q}\ldots s_{n}  } \ q_{s_{i_1}} q_{s_{i_2}}
\ldots q_{s_{i_q}} \ \prod _{j\in \{1, \ldots, n\}- \{i_1, i_2,
\ldots, i_q\}}
 q_{s_j}^{2k_j} $$
o\`u les $k_j$ sont des entiers naturels.}

\bigskip  La preuve se fait pas r\'ecurrence sur $n$. Lorsque
$n=1$ c'est
\'evident. Supposons cette propri\'et\'e vraie pour $n-1$ et montrons
la pour $n$. Si $i_1>1$, la suite 
$(s_{i_1},s_{i_2},\ldots, s_{i_q})$ est  extraite 
de la suite $(s_2,\ldots, s_n)$ et elle est $ws_1$-distingu\'ee.
L'assertion   est vraie pour
 $(ws_1, s_2,\ldots ,s_n)$ par hypoth\`ese de r\'ecurrence, et l'on
obtient le r\'esultat, avec $k_1=0$. Supposons que $i_1=1$.
  La suite 
$( s_{i_2},\ldots, s_{i_q})$ est  extraite 
de la suite $(s_2,\ldots, s_n)$ et elle est $w$-distingu\'ee.
L'assertion est vraie pour $(w, s_2,\ldots , s_n)$ par hypoth\`ese de
r\'ecurrence, et l'on obtient

$$q_{ws_2\ldots s_n}\ =  \ 
q_{w s_2\ldots   \hat s_{i_2} \ldots
\hat s_{i_q}\ldots s_{n}  } \  q_{s_{i_2}}
\ldots q_{s_{i_q}} \ \prod _{j\in \{2, \ldots, n\}- \{ i_2,
\ldots, i_q\}}
 q_{s_j}^{2t_j} $$
pour des entiers $t_j\geq 0$. L'assertion pour $i_1=1$  
r\'esulte alors du lemme suivant.

  \bigskip   {\bf 1.3 \ Lemme \ } {\sl  Soit $w   \in W$  et
$s_1, s_2, \ldots s_n \in S$ tels que $w<ws_1$ et $s_1\ldots s_n$
est une d\'ecomposition r\'eduite. Alors

a) $ws_2\ldots s_n < ws_1s_2\ldots s_n$

b) $$q_{ws_1\ldots s_n}   q_{ws_2\ldots s_n }^{-1}\ =  \ 
q_{s_1} \   \prod _{j\in \{2, \ldots, n\} }
 q_{s_j}^{2k_j} $$
pour des entiers $0\leq k_j \leq 1$.}

\bigskip Le lemme 1.3  termine la d\'emonstration  du
lemme fondamental lorsque
$v\in W_{aff}$. Le cas g\'en\'eral s'en d\'eduit, car si
$v=v_{aff}u, v_{aff}\in W_{aff}, u\in \Omega$, et si $w,w'\in W$,
on a $ \tilde T_w \tilde T_{v^{-1}}^{-1}= \tilde T_w \tilde
T_{v_{aff}^{-1}}^{-1}T_u, \ 
 T_{w'} T_{u}=T_{w'u}, \ q_{w'u}=q_{w'},  \ w'<wv_{aff} $ implique
$w'u<wv$. 
La d\'emonstration du lemme fondamental sera donc achev\'ee, apr\`es
la preuve de 1.3.

\bigskip Preuve du lemme 1.3. La partie a) est d\'emontr\'ee dans 
([Ha] lemma 5.6). Elle implique que le membre de gauche de b) est un  
mon\^ome
$c_n\in Q$ diff\'erent de $1$. On montre b) par r\'ecurrence sur $n$.
Pour $n=1$ on a $c_1=q_{s_1}$ donc b) est vrai. 
Montrons que si b) est vrai pour $n-1$, alors il est vrai pour $n$. On
a 

$q_{ws_1\ldots s_n}=  q_{ws_1\ldots s_{n-1}}q_{s_n}^{\e_1}$

$q_{ws_2\ldots s_n}=  q_{ws_2\ldots s_{n-1}}q_{s_n}^{\e_2}$

\noindent avec $\e_1,\e_2\in \{ \pm 1 \}$. On a donc 

$c_n\ =
  \ c_{n-1} , \ c_{n-1} q_{s_n}^2, \ c_{n-1}
q_{s_n}^{-2}$ selon que  $\e_1=\e_2, \ (\e_1,\e_2)=(1,-1), \
(\e_1,\e_2)=(-1,1)$.

\noindent Par hypoth\`ese de r\'ecurrence, 

$c_{n-1}=q_{s_1} \   \prod _{j\in \{2, \ldots, n-1\} }
 q_{s_j}^{2t_j} $ 
pour des entiers $0\leq t_j \leq 1$. 

Donc $c_n$ v\'erifie  b) si $c_n= c_{n-1}$ ou $   c_{n-1} q_{s_n}^2$.
Pour terminer la d\'emonstration il suffit de remarquer que si
$c_nq_{s_n}^2\ =c_{n-1}$ alors il existe $j$ tel que $t_j=1$ et
$q_{s_n}=q_{s_j}$. Donc $c_n$ v\'erifie aussi  b). 

Le lemme 1.3 est d\'emontr\'e.

\bigskip  {\bf D\'ecomposition de Bernstein \ }

 D'apr\`es Bernstein-Lusztig, la  d\'ecomposition de $W$ en
produit semi-direct 
$W=W_o e^X$ se g\'en\'eralise  \`a
l'alg\`ebre de Hecke
 $H_{\phi_L}[q_*^{-1/2}]$. 

On notera
$?_x:=?_{e^x}$  pour
$x\in X$, et $A$
l'image du morphisme de
$\Z[q_*^{\pm 1/2}] $-modules 
$$\tilde
\theta:\Z[q_*^{\pm 1/2}][X]\to H[q_*^{- 1/12}]$$ qui envoie $x\in
X$ sur
 $\tilde \theta_x:= \tilde T_y\tilde
T_z^{-1}=q_y^{-1/2}q_z^{1/2}T_yT_z^{-1}$   o\`u
  $x =y-z$ pour deux \'el\'ements dominants $y,z$ de $X$ (bien
d\'efini car la longueur
\'etant un morphisme de monoide sur les \'el\'ements dominants).

\bigskip Soit  $\alpha
\in B$ une racine simple de r\'eflexion associ\'ee $s\in S_o$. Si 
$\alpha^\vee  \in 2X^\vee$, soit $\tilde s \in S_o$ la 
r\'eflexion d\'efinie  comme en [Lusztig1 2.4 page 604]. On d\'efinit
deux \'el\'ements $a_s, b_s
\in A$  par:

$$a_s=q_s-1, \ \ b_s= 1-\tilde \theta_{-\alpha}, \ \ {\rm si \
} \alpha^\vee \not\in 2X^\vee,$$

$$a_s=q_s-1 +q_s^{1/2}
( q_{\tilde s  } ^{1/2}- q_{\tilde s  }^{-1/2})\tilde \theta_{-\alpha}
,
\
\ b_s= 1-\tilde
\theta_{-2\alpha}, \ \ {\rm si
\ } \alpha^\vee  \in 2X^\vee.$$
 
\bigskip {\bf 1.4  \ Th\'eor\`eme \ } (Bernstein-Lusztig) {\sl

1)  Le morphisme $\tilde \theta:\Z[q_*^{\pm 1/2}][X]\to H[q_*^{-
1/12}]$ est injectif et respecte le produit.

2) \ $ \tilde
\theta_x  T_s - T_s  \tilde \theta_{s(x)}= ( \tilde
\theta_x - \tilde \theta_{s(x)})a_s b_s^{-1}
 $ pour tout $(s,x)\in S_o\times X$.

3) $(T_{w_o})_{w_o\in W_o}$ est une base de  
$H[q_*^{- 1/2}]$ comme $ A$-module, \`a droite ou \`a gauche.  

4)    Le
centre de
$H[q_*^{- 1/2}]$ est l'alg\`ebre $A^{W_o}$ des \'el\'ements
$W_o$-invariants de
$ A$. C'est un $\Z[q_*^{\pm
1/2}] $-module libre de base $\tilde z_x=\sum _{x'\in W_oe^x}\tilde
\theta_{x'}$ pour tout $x\in X_{dom}$.

}

 \bigskip 
Le membre de gauche de la partie 2) du th\'eor\`eme est un
\'el\'ement de $A$, dont le produit avec $b_s$
est
$(
\tilde
\theta_x - \tilde \theta_{s(x)})a_s$. 

\bigskip  Les r\'ef\'erences que je connais  concernent seulement la
$\Z[v^{\pm 1}]$-alg\`ebre $H_{\phi_L}[v^{-1}]$.  Ce sont [Lusztig1 3.3
\`a 3.11 page 607-610], [Lusztig2 7.1 page 216, 8.1 page 222].
 L'alg\`ebre $H_{\phi_L}[v^{-1}]$ suffit pour  les applications en vue,
et  il me semble que les d\'emonstrations restent valables pour
$H[q_*^{-1/2}]$, modulo des modifications triviales.

\bigskip Par le th\'eor\`eme,  $H[q_*^{- 1/2}]$ est un $  A$-module
\`a droite ou
\`a gauche libre de rang fini \'egal au cardinal $|W_o|$ de $W_o$.
{\sl Si de plus, $X$ contient tous les poids dominants pour le
syst\`eme de racines $R$, alors $\Z[X]$ est libre sur
$\Z[X]^{W_o} $ de rang
$|W_o|$ {\rm [Steinberg th.2.2 page 173]} et
$H[q_*^{- 1/ 2}]$ est un module libre de rang $|W_o|^2$ sur son
centre.}

\bigskip A-t-on des r\'esultats analogues pour la
$\Z[q_*]$-alg\`ebre $H$ ?  
\bigskip 
- Le $\Z[q_*]$-module
$A\cap H$ est-il libre ? 

- L'alg\`ebre commutative $A\cap H$
  est-elle de type fini ? 

- Le centre $Z(H)=A^{W_o}\cap H$ de $H$ 
est-il une alg\`ebre de type fini ?

- Le $Z(H)$-module $A\cap H$ est-il  
 de type fini ?

-  Le $(A\cap H)$-module $H$ {\bf \`a droite}, est-il
 de type fini ?

\bigskip  Toutes ces questions seront r\'esolues avec l'aide du
th\'eor\`eme suivant qui se d\'eduit du lemme fondamental appliqu\'e
\`a
 $\tilde T_{w_o} \tilde \theta_{y-z}=\tilde T_{w_o e^y} \tilde
T_{e^z}^{-1}
$ pour tout 
$w_o\in W_o, \ y,z\in X_{dom} $ (l'\'egalit\'e vient de
$\ell(w_oe^y)=\ell(w_o)+\ell(e^y)$).

\bigskip {\bf 1.5  \ Th\'eor\`eme }\ {\sl Pour  $ 
w_o\in W_o, \ x\in X $, posons

$$E_ {w_oe^x}:= q_{w_oe^x}^{1/2} \tilde T_{w_o}\tilde \theta_x. $$ Pour
tout
$w\in W$, le d\'eveloppement de $E_w$ dans la base de Iwahori-Matsumoto
$(T_{w })_{w\in W}$ est
 
$$E_w=T_w+\sum_{ w'<w} \ a_{w'} \ T_{w'}  
$$ 
  pour des polyn\^omes $a_{ w'} \in \Z[q_*]$. }

 \bigskip 
Le th\'eor\`eme 1 s'en d\'eduit  par un 
argument g\'en\'eral.

\bigskip Le th\'eor\`eme 2 r\'esulte des des propri\'et\'es
remarquables ci-dessous de la nouvelle base $(E_w)_{w\in W}$ de
l'alg\`ebre de Hecke g\'en\'erique $H$.

\bigskip {\bf 1.6  \ Propri\'et\'es de  $(E_w)_{w\in W}$ \ }

\bigskip {\bf (1.6.1) }  
{\sl  Le $\Z[q_*]$-module libre
de base $(E_x)_{x\in X}$ est \'egal \`a la sous-alg\`ebre commutative
$A\cap H$.} 

C'est imm\'ediat puisque $E_x=q_x^{1/2} \tilde \theta_x$.
\bigskip 
{\bf (1.6.2) }    
{\sl On a  $H=\oplus_{w_o\in W_o}
H(w_o) $ 
o\`u $H(w_o)$ est le $\Z[q_*]$-module
libre
 de base $(E_{w_oe^x})_{x\in X}$.} 

\bigskip En effet, 
$\tilde T_{w_o}\tilde \theta_x\tilde \theta_{x'}=\tilde
T_{w_o}\tilde \theta_ {x+x'}$ pour tout
$ x,x'\in X$, car $\tilde \theta$ respecte le produit
(1.4.1).

\bigskip 

{\bf (1.6.3) } {\sl  $H(w_o)=T_{w_o} J(w_o),$  o\`u 
$J(w_o)$ est un id\'eal fractionnaire de type fini de
 $A\cap H$.}

\bigskip  En effet,
$$E_{w_oe^x }E_{x' }=  q(w_oe^x, e^{x'})
 E_{w_oe^{x+x' }}, \ \ \  q(w_oe^x, e^{x'}):=(q_{w_oe^x}
q_{ x' }q_{w_oe^{x+x'}}^{-1})^{1/2}.  
 $$
Le produit $q(w_oe^x, e^{x'})\in \Z[q_*]$ est un mon\^ome  
puisque $(E_w)_{w\in W}$ est une base de $H$ sur $\Z[q_*]$.  
On montre 
avec la formule des longueurs: 

{\sl Il  existe un 
 ensemble fini $X(w_o)=\{x_1, \ldots, x_r\} \subset X$ d\'ependant de
$w_o$,  tel que  
$$X=\cup_{i=1}^r X(w_o,x_i)$$ o\`u   $X(w_o,x_i)$ est
l'ensemble des
$x\in X$ v\'erifiant
 $\ell(w_oe^x)=\ell(w_oe^{x_i})+\ell(e^{x-x_i})$.} 

\bigskip

Si $ x\in  
X(w_o,x_i)$, on a 
$$E_{w_oe^x}=E_{w_oe^{x_i}}E_{x-x_i},  \ \ \ q(w_o,e^{x_i})
E_{w_oe^{x_i}}= T_{w_o}E_{x _i},  \ \ \ E_{w_oe^x}=c (w_o,e^x)
T_{w_o}E_x,$$ o\`u
$q(w_o,e^{x_i})\in \Z[q_*]$ est un mon\^ome   et 
$ c
(w_o,e^x)=q(e^{x_i}, e^{x-x_i})q(w_o,e^{x_i})^{-1}  $ 
 est un mon\^ome de Laurent. On a $ q(w_o,e^{x_i})E_{w_oe^{x_i} }=
T_{w_o} E_{x_i}
 $.

{\sl -  Le $A\cap
H$-module \`a droite
$H(w_o)$ est engendr\'e par
  $(E_{w_oe^{x } })_{x\in X(w_o)}$, 

-    $H(w_o)=T_{w_o}J(w_o)$
o\`u $J(w_o)$ est le
$A\cap H$-module   de type fini engendr\'e par
$(q(w_o,e^{x})^{-1}E_{x})_{x\in X(w_o)}$ dans
$A\cap H(q_*^{-1})$. }

\bigskip {\bf (1.6.4) } On a $E_xE_{x'} =E_{x+x'}$ s'il existe $w_o\in
W_o$ tel que
$w_o(x),w_o(x') \in X_{dom}$ car la longueur est $W_o$-invariante sur
$e^X$ et un morphisme de monoide sur $e^{X_{dom}}$. Le
monoide 
$X_{dom}$ est  de type fini. On choisit un syst\`eme g\'en\'erateur
fini
$M_{dom}$ de 
$X_{dom}$. Alors  
 $M_X:=W_o(M_{dom})$ est un ensemble fini.

{\sl La
$\Z[q_*]$-alg\`ebre 
$A\cap H$ est de type fini  engendr\'ee par $(E_x)_{x\in
M_X} $.}

\bigskip {\bf (1.6.5) } Le centre $Z(H)$ de $H$ est $A^{W_o}\cap
H$. Comme 
$q_x=q_{w_o(x) }$  pour tout
$(w_o,x)\in W_o\times X $, on a
$q_x^{1/2}\tilde z_x=\sum_{x'\in W_o(x)} E_{x'} $. On
d\'eduit alors de (1.4.3):

 {\sl Le centre $Z(H)$ de $H$ est \'egal \`a $ (A\cap
H)^{W_o}$; c'est un $\Z[q_*]$-module libre  de base
$ Z_x:=\sum_{x'\in W_o(x)} E_{x'}$ pour $x\in X_{dom}$.
  }

\bigskip {\bf (1.6.6) }    Des arguments g\'en\'eraux
 [Bourbaki AC V 1.9 th\'eor\`eme 2 page 29]  permettent
de d\'eduire:

 {\sl $A\cap H$ est  un 
$(A\cap H)^{W_o}$-module de type fini}, car c'est une
$\Z[q_*]$-alg\`ebre de type fini  (1.6.4).  

{\sl $(A\cap H)^{W_o}$
  est une $\Z[q_*]$-alg\`ebre de type fini}, car 
$\Z[q_*]$ est un anneau noetherien.

\vfill\eject \centerline {\bf Chapitre 2 \ }

\bigskip \centerline{{\bf Caract\`eres de $A\cap H$ et modules
standards}}

\bigskip   Soit $R$ un anneau commutatif int\`egre de corps des
fractions
$K$ 
 et $i :
R\to
K$ l'inclusion.   
Soient    un
morphisme d'anneau de $\phi:\Z[q_*  ]\to R $ et   $\X:A\cap H \to R$ 
un morphisme d'anneau     prolongeant
$\phi$. On note $\w$  la 
restriction  de $\X$ au centre $Z(H)$  de $H$. 

\bigskip Pour les applications
que nous avons en  vue, il ne serait pas g\^enant de supposer que
$\phi$ se prolonge \`a
$ \Z[q_* ^{1/2} ]$.  On notera $B_{\phi}:=B\otimes_{\phi}R$
la
$\phi$-sp\'ecialisation d'une
$\Z[q_*]$-alg\`ebre $B$ et $f_{\phi}$ la
$\phi$-sp\'ecialisation d'un morphisme $f$ de $\Z[q_*]$-alg\`ebre.

 \bigskip On identifie $\X$ 
\`a un morphisme  d'anneau 
$\X_{\phi}:(A\cap H)_{\phi}
\to R$, et aussi en posant $ \X_X(x):=\X(E_x)$, \`a une application 
$\X_X: X\to R $ v\'erifiant 
$$\X_X(x) \X_X(x') = \phi
((q_{x}q_{x'}q_{x+x'}^{-1})^{1/2})\X_X(x +x') 
$$ pour tout $x,x'\in X ,$ par la formule
$E_xE_{x'}=(q_{x}q_{x'}q_{x+x'}^{-1})^{1/2} E_{x+x'}$
cas particulier de (1.6.2).  

\bigskip Lorsque $\phi(q_s)\subset R^*$ est inversible pour tout $s\in
S$, et que
$\phi$ se prolonge
\`a
$ \Z[q_*^{1/2}]$, en posant pour $x\in X$, $$ \X_X(x)=  \X_A
(\tilde\theta_x):=
\phi(q_x)^{-1/2}\X(E_x),$$ on identifie $\X$ \`a un morphisme
d'anneau
$ \X_A: A\to R$  et \`a 
 un  morphisme de groupe
 $\X_X:X \to R^*$.

 \bigskip  On dit que $\w$ est le
caract\`ere central d'un
$H$-module $V$, lorsque $Z(H)$ agit sur $V$ par
multiplication par $\w$. Le caract\`ere $\w$ s'identifie  
\`a un caract\`ere de $Z(H)_{\phi}$. Le centre de la
$\phi$-sp\'ecialisation $H_{\phi}$ de $H$ peut
\^etre plus gros que  $Z(H)_{\phi}$. 
La $R$-alg\`ebre  $H_{\phi}$ est un $R$-module libre  de rang infini
puisque $W$ est infini. La 
$R$-alg\`ebre
$H_{\w}=H\otimes_{\w} R$  est  de type fini sur $R$,
 car $H$ est un module de rang fini sur son centre. On en d\'eduit:

\bigskip {\bf 2.1 \ Proposition \  } {\sl Si $R$ est un corps
alg\'ebriquement clos,  un
$H_{\phi}$-module  simple  est   de dimension finie si et
seulement s'il a un caract\`ere central}.

\bigskip 
 {\bf\ 2.2  \ D\'efinition \ }{\sl   Le
$H_{\phi}$-module  \`a gauche induit de $\X$
$$I(\X ):= H \otimes_{  \X } R $$
de caract\`ere central $\w$ est 
 appel\'e un
  module standard \`a gauche de $H_{\phi}$. }

\bigskip Lorsque $\phi(q_s)\subset R^*$ est inversible  pour tout $s\in
S$,  et que
$\phi$ se prolonge
\`a
$ \Z[q_*^{1/2}]$, le module standard $I(\X)$ est le module
standard  $I(\X_A)$ induit du caract\`ere
$\X_A:A\to R$ prolongeant $\X$.
Il est libre sur $R$ de rang $|W_{o}|$ et admet une base canonique:   l'image canonique de $(T_{w_{o}})_{w_{o}\in W_{o}}$.

\bigskip {\bf 2.3  \ Premi\`eres propri\'et\'es \ }  Le  module
standard
$I(\X)$ satisfait les propri\'et\'es  suivantes.

\bigskip   {\bf  (2.3.1)  \ }  {\sl Il
est de type fini  comme  
$R$-module,} car $H$ est de type fini comme $A\cap
H$-module (th\'eor\`eme 2).

\bigskip {\bf  (2.3.2)  \ } {\sl Il
est   
 cyclique, avec un g\'en\'erateur canonique $1\otimes 1$ propre pour
$A\cap H$ de valeur propre $\X$. Un $H_{\phi}$-module 
\`a gauche  engendr\'e  par un vecteur   propre pour $A\cap H$ de
valeur propre $\X$ est quotient  de $I(\X)$.
 }

  \bigskip Le th\'eor\`eme 3 se d\'eduit de ce r\'esultat et de (2.1).   Comme
 la $R$-alg\`ebre de type
fini  $H_{\w}$  admet un module simple \`a gauche, on obtient aussi:
 
\bigskip {\bf (2.3.3) \  } {\sl Si $R$ est un corps
alg\'ebriquement clos,  tout morphisme de
$R$-alg\`ebre  $\w:Z(H ) \to R$ se prolonge en un morphisme d'anneau 
$\X:A\cap H\to R$.}

\bigskip  Supposons   $\X(q_{s})\neq 0$ pour tout $s\in S$.
Par d\'efinition, une structure enti\`ere du module standard $I(\X,K)=I(\X)\otimes_{R }K$ est un sous-$R$-module de type fini stable par $H$   et contenant une base de $I(\X,K)$. Le noyau  du morphisme  canonique   $ I(\X)\to I(\X,K)$ est  le sous-groupe de torsion $I(\X)_{tor}$.

\bigskip {\bf (2.3.4) \  } {\sl L'image de $I(\X)$ dans $I(\X,K)$ est une structure enti\`ere $M(\X)\simeq I(\X)/I(\X)_{tor}$  contenant la base canonique.}

\bigskip Le morphisme $I(\X)\to I(\X,K)$ est le produit tensoriel avec $H$ de l'inclusion  $i: R \to K$, vue comme une inclusion  de $A\cap H$-modules via le caract\`ere $\X$. Si $H$ est un $A\cap H$-module \`a droite plat, alors le morphisme  est injectif et $I(\X)$ n'a pas de torsion.
La classification des $R$-modules de type fini  implique que si $R$ est un anneau principal, 
 le rang du $R$-module 
$I(\X)$ est $|W_{o}|$ et le $R$-module $M(\X)$ est libre de rang $|W_{o}|$.

\bigskip {\bf 2.4 \ Proposition \ } {\sl  Supposons   $\X(q_{s})\neq 0$ pour tout $s\in S$ et que
$R$ est un anneau de valuation discr\`ete de corps r\'esiduel $k_{R}$, et  de r\'eduction $r:R\to k_{R}$. Alors  les propri\'et\'es suivantes sont \'equivalentes:
 
 a) Le module standard $I(\X,k_{R}):= H\otimes_{r \X}k_{R}=I(\X)\otimes_{r} k_{R}$ est de dimension $|W_{o}|$ sur $k_{R}$,
 
 b) La r\'eduction de la structure enti\`ere canonique de $I(\X,K)$ est le 
 le module standard $I(\X,k_{R})\simeq  M(\X)\otimes_{R,r}k_{R}$,
 
 c) Le $R$-module $I(\X)$ est sans torsion.}
 
\bigskip Le th\'eor\`eme 5 r\'esulte de cette proposition, car le cas $R=\overline \Z_{p}$ se ram\`ene au cas de l'anneau des entiers d'une extension finie de $\Q_{p}$.

\bigskip {\bf 2.5 \ Remarques  \ }  (i) Avec les hypoth\`eses de la proposition, il  existe une application $f:W_{o} \to X$, non unique et d\'ependant de $\X$, telle que  $M(\X)$ 
admet comme base sur $R$ l'image canonique de $(E_{w_{o}e^{x_{o}} })_{w_{o}\in W_{w_{o}}}$.
Lorsque  le $R$-module $I(\X)$ est sans torsion (par exemple si $r \phi (q_{s})\neq 0$ pour tout $s\in S$), l'on peut probablement montrer par r\'eduction, en utilisant sur $R$ les arguments classiques [Rogawski], [Cherednik, proposition 1.5 page 416],  que pour tout $w_{o}\in W_{o}$, 

- $I(\X, k_{R})$ et $I(\X w_{o}, k_{R})$ ont les m\^emes suites de Jordan-H\"older,

- $I(\X, k_{R})$ contient un vecteur propre pour $A\cap H$ de valeur propre $\X w_{o}$, ou ce qui est \'equivalent, il existe un $R$-morphisme $H$-\'equivariant (un op\'erateur d'entrelacement)
 de  $I(\X w_{o}, k_{R})$ dans $I(\X, k_{R})$.
 
(ii) Supposons que $R=\overline \Z_{p}$ et que la r\'eduction de $\phi(q_{s})$ est nulle pour tout $s\in S$.
Que peut-on dire sur la  r\'eduction du  sous-module de torsion $I(\X)_{tor}$ ? 
 Est-ce que tout caract\`ere $A\cap H \to \overline \F_{p}$ nul sur  $q_{s}$ pour tout $s\in S$, se rel\`eve ?

\vfill\eject {\centerline{\bf Chapitre 3 \ Exemple du groupe
lin\'eaire}}

\bigskip  Soient $n>1$ un entier et
$q$ une ind\'etermin\'ee.  On note $A_n$  la
$\Z[q]$-alg\`ebre commutative engendr\'ee par les \'el\'ements
$(e_I)_{ 
 I\subset \{1, \ldots, n \}}$ v\'erifiant les
relations  
$$e_\emptyset =1, \ \ e_Ie_J=q^{yz}e_{I \cup J}e_{I \cap J}, \ \ \
|I\cap J|=x, |I |=x+y, \ |J|=x+z,$$ 
 en notant $|I|$ le nombre d'\'el\'ements de
$I$. 
L'identit\'e polyn\^omiale  
$$(X+Y)(X+Y-1)  + (X+Z)(X+Z-1)+2YZ  = (X+Y+Z)(X+Y+Z-1)  + X(X-1) $$
permet de v\'erifier que {\sl les relations sont \'equivalentes \`a  
 $$e_\emptyset =1, \ \ e_{i_1}e_{i_2}\ldots e_{i_t}=q^{t(t-1)/2}
e_{\{i_1,
\ldots, i_t\}}. $$
 }
 Le groupe $S_n$ des permutations  de  $\{1, \ldots, n \}$, permute les
parties de $\{1, \ldots, n \}$ de m\^eme cardinal, et agit
naturellement sur 
$A_n$ en fixant $e_{\{1, \ldots, n\}}$.  L'alg\`ebre $A_n^{s_n}$ fixe
par $S_n$  est engendr\'ee  par $ z_t:=\sum _{|I|=t}e_I$
pour $1\leq t \leq n $.

 On pose  $z :=z_n=e_{\{1, \ldots, n
\}}$.

\bigskip {\bf 3.1 \ Proposition \ }  {\sl Pour la donn\'ee radicielle
bas\'ee correspondant au groupe $GL(n)$, l'alg\`ebre commutative
$A\cap H$  est  isomorphe \`a  
$A_n[ z^{-1}] $.
Le centre   de
$H$  est isomorphe \`a $A_n^{S_n}[ z^{-1}]$.}
\bigskip Preuve.  La donn\'ee radicielle se d\'ecrit avec le
   groupe
$GL(n,\C)$. Soit $e\in \C$ un \'el\'ement transcendant.
 Le groupe
ab\'elien $ X=\Z^n$ s'identifie au groupe des 
diagonales
$e^x:=\diag (e ^{x_1},
\ldots, e ^{x_n})$ pour  $x=(x_1, \ldots, x_n)\in X$. La base
$B=(\alpha_i)_{1\leq i \leq n-1}$   est telle que
  $$e^{\alpha_i}=(1,\ldots,1,e^{-1},e,1, \ldots, 1),  \  \ \
\alpha_i^\vee
\ (x)\ =\ x_{i+1}-x_i$$ pour
$1\leq i
\leq n-1$, o\`u $e^{-1}$ est \`a la place $i$.  Le groupe de Weyl fini
$W_o$ s'identifie au groupe sym\'etrique $S_n$ engendr\'e par les
r\'eflexions  
$s_i=(i,i+1)$ pour $1\leq i \leq n-1$ associ\'ees aux racines
$\alpha_i$. Le groupe de Weyl  
$W$ de la donn\'ee radicielle est isomorphe au produit semi-direct 
$ S_n. e^X$. Le sous-groupe
$\Omega\subset W$ des
\'el\'ements de 
 longueur $0$  s'identifie au  groupe $t^{\Z}$  o\`u 
  $$t=\pmatrix {0&1&0&0&\ldots &0\cr 0&0&1&0&\ldots &0 
\cr   & & & & &  \cr
\ldots &\ldots &\ldots &  \ldots &\ldots &\ldots \cr      & & & & &
\cr  0&1&0&0&\ldots &1\cr e&0&0&0&\ldots &0&\cr}. $$
 Le groupe de Weyl affine $W_{aff}$ est engendr\'e par les
r\'eflexions
$s_o, s_1, \ldots , s_{n-1}$, avec $s_ot=ts_1, s_1t=ts_2$,
 $\ldots, s_{n-2}t=ts_{n-1}$. 
Les  poids fondamentaux dans  $ X $ sont $\w_i
= (0^i, 1^{n-i})$,  
$i$ premiers termes
\'egaux \`a $0$ et $n-i$ derniers termes \'egaux \`a $1$, 
pour $ 1\leq i \leq n-1$.  Il est
pratique d'introduire   $
\w_n:=(0^n), \ \w_o:=(1^n)$. On a $X\cap \Omega =\w_o \Z $ et  le
centre de $GL(n,\Q_p)$ est
$p^{\w_o\Z}$. La
longueur de $\w_i$ est  le nombre  $(n-i)i$ de racines positives
contenant
$\alpha_i$ si $1\leq i \leq n-1$. Elle est nulle si $i=0,n$.  
  Les  conjugu\'es de $\w_{n-1} $ par $ 
S_n$  sont  
$y_i:= \w_{i-1}- \w_i 
\ (1\leq i
\leq n)$,   avec un seul coefficient $1$ diff\'erent de $0$  
\`a la place $i$. 

\bigskip  Une
base du monoide $X^+$ est $(w_i)_{1\leq i \leq n}$. L'alg\`ebre
$A\cap H$ est donc engendr\'ee par $E_{w_o(w_i)}$ pour $w_o\in S_n, \
1\leq i \leq n$. On a
$y_n=\w_{n-1}
$ et
$ y_{j}=s_{j}\ldots  s_{n-1}(\w_{n -1})$ pour $1\leq j \leq n-1$.
 Plus g\'en\'eralement, les
 $S_n$-conjugu\'es de $\w_{n-t}$ sont $y_I:=\sum_{i\in I}y_i$ pour
toutes les parties $I$ \`a
$ t$ \'el\'ements  de
$\{1,\ldots,n\}$.  On note pour simplifier $E_I=E_{y_I}, \tilde \theta_I=\tilde
\theta_{y_I}$. L'alg\`ebre
$A\cap H$ est   engendr\'ee par $E_I$ pour les parties non vides de
$\{1,\ldots,n\}$. 

Les relations entre les $E_I$ se d\'eduisent de celles entre les
$\tilde \theta_I$.  On a   $q_{\w_{n-t}}=q_{I}=q^{t(n-t) }$
 pour
$I\subset
\{1,
\ldots, n\}$ avec $t$ elements, ainsi:
$$\tilde \theta_{i}= q^{-i + (n+1)/2} T_{\w_{i-1}}T_{\w_i}^{-1}, \
\
\  E_i=q^{n-i} T_{\w_{i-1}}T_{\w_i}^{-1}
\ \ , \ \ \ 
\tilde \theta_I =\prod_{i\in I} \tilde \theta_i, \
\ E_I= q^{ t(n-t)/2}
\tilde \theta_I. $$ 
  L'\'el\'ement central inversible $E_{\{1,
\ldots, n\}}=T_{\w_o}
$ est not\'e plus simplement
$Z$.   On a  
$$E_\emptyset = 1,\ 
\
\prod_{i\in I}E_{i} = q^{t(t-1)/2}E_I. $$
 La $\Z[q_*]$-alg\`ebre engendr\'ee par    $Z_t:=\sum_{|I|=t} E_I\
(1\leq t
\leq n-1), Z^{ \pm 1}$ est \'egale  au centre $Z(H)$.

On  en d\'eduit la proposition 3.1.

\bigskip Soit $R$ un anneau int\`egre.
Un 
morphisme d'anneau $ \X :A\cap H \to R$   v\'erifie
$$\prod_{i\in I} \X
(E_i) =\X(q)^{t (t-1)/2}  \X
(E_I),$$    
 pour tout $I\subset \{1,\ldots,
n\}$ avec $t\geq 1$ \'el\'ements. 
Si $\X(q)$ est inversible
dans
$R$, alors
$  \X$ est d\'etermin\'e par les $n$ \'el\'ements non nuls $ 
(\X(E_i))_{1\leq i
\leq n }$ de produit $\prod_{1\leq i \leq n
}  \X(E_i)$   inversible, et inversement.  Lorsque $\X(q)=0$, la
situation est bien diff\'erente.

\bigskip {\bf 3.2  \ Lemme \ } {\sl Soit $\X :A\cap H \to R$ un
morphisme  d'anneau tel que $\X(q)=0$.  

1) Les parties non vides 
$I\subset \{1, \ldots, k\}$ telles que $\X(E_I) \neq 0$ forment une
suite  croissante 
$$\emptyset \neq I_1\subset I_2 \subset \ldots \subset
I_{r}=\{1,
\ldots, n\}.
$$

2) Inversement, \'etant donn\'e  une suite croissante comme ci-dessus 
et des
\'el\'ements non nuls $x_{I_1}, \ldots, x_{I_r}
\in R$, il existe un unique morphisme
d'anneau $\X:A\cap H\to R$ tel que $\X(E_{I_j})=x_{I_j}$ et
$\X(E_I)=0$ si $I\neq I_j$ pour tout $1\leq j \leq r$.}

\bigskip  Preuve.   On choisit $x_I\in R$ pour toute partie non vide
$I$  de
$\{1,
\ldots, n \}$. Il existe un caract\`ere $\X :A\cap H \to R$  tel que
$\X(q)=0$ et $\X(E_I)=x_I$ pour tout $I$, si et seulement si
$x_{\{1\ldots, n\}}\in R^*$ est une unit\'e et $x_Ix_J=0$ lorsque
$yz\neq 0$ o\`u 
$|I|=x+y,|J|=x+z,|I\cap J|=x$, i.e. lorsque $I$ n'est pas contenu dans
$J$  et
$J$ n'est pas contenu dans $I$.

Deux parties $I,J$ tels que $x_Ix_J\neq 0$ doivent donc v\'erifier
$I\subset J$ ou $J\subset I$. Les parties $I$
non vides telles que $x_I\neq 0$ doivent former une suite croissante,
se terminant par $\{1\ldots, n\}$. Il n'y a aucune autre condition sur
les valeurs non nulles
$x_I$ ni sur la suite croissante.

Le morphisme $\X$ est \'evidemment
d\'etermin\'e par les $x_I$, lorsqu'il existe. 

On en d\'eduit le lemme 3.2.

\bigskip {\bf Applications \ } 1)  Le nombre de $S_n$-conjugu\'es de
$\X
$   est 
$${  n!\over t_1!(t_2-t_1)! \ldots (n-t_{r-1}) !} $$
 si  $t_j$ est le nombre d'\'el\'ements de  $ I_j $ pour $1\leq j
\leq r$.

 2) La restriction de $\X$ au centre $Z(H)$ est le caract\`ere $\w$
tel que   
  $$ \w (Z_{t_j})=\X (E_{I_j})   \ \ \ {\rm pour \ tout } \ \ \ 
1\leq j
\leq r ,$$ et
$\w(Z_t)=0$ pour les autres valeurs possibles de $t$.

\bigskip {\bf 3.3 \  Definition \ }{\sl Soit $R$ un anneau int\`egre.
On dit qu'un $R$-caract\`ere $\X$ de $A\cap H $ nul sur $q$, ou que sa
restriction 
$\w$ au centre $Z(H)$,   est

- r\'egulier,  si  le
stabilisateur de $\X$ dans $S_n$ est trivial, i.e. $\w(Z_t)\neq 0$
pour tout $1\leq t \leq n$,

- singulier,   s'il n'est pas r\'egulier,

- supersingulier,    si  $\X$ est fixe
par $S_n$, i.e.
$\X (E_I)=0$  pour tout  $\emptyset \neq I\subset \{1,\ldots,
n\}$, i.e. $\w(Z_t)=0$ sauf si $t=n$. 

  Un $H_{r_p\phi}$-module ayant un caract\`ere central
 est appel\'e r\'egulier, ou singulier, ou supersingulier, si son
caract\`ere central a cette propri\'et\'e.} 

\bigskip  La terminologie
``supersingulier''est analogue
\`a  celle introduite par   Barthel et Livne dans leur
\'etude des
$\overline
\F_p$-repr\'esentations de
$GL(2,F)$ sur un corps
$p$-adique $F$. 
 Un caract\`ere 
$\w,
\X$ supersingulier est   d\'etermin\'e  par sa valeur 
$z:=\X (Z )=\w(Z)
\in
R^*$ en $Z$, et sera not\'e  $\w_z, \X_z$.

\bigskip {\bf 3.4  \  Preuve du th\'eor\`eme 6 \ }  Soit $\phi(q
)\in \overline \Z_p$ non nul de r\'eduction
nulle $r_p\phi(q)=0$.   D\'emontrons
  que tout morphisme d'anneau
$\X :A\cap H
\to
\overline
\F_p$ tel que $\X(q)=0$ se rel\`eve en un morphisme d'anneau $\tilde \X
:A\cap H
\to
\overline
\Z_p$ tel que $\tilde \X (q) =\phi(q)$.

\bigskip    On choisit un syst\`eme compatible de racines
$n$-i\`emes
$\phi(q)^{1/n} $ de $\phi(q)$ dans $\overline \Z_p$, pour tout entier
$n\geq 1$. Tout
\'el\'ement $\lambda_i$ de
$\overline
\Z_p$ s'\'ecrit de facon unique $\lambda_i= x_i\phi(q)^{e_i}$ pour un
nombre rationnel
$e_i\geq 0$ et une unit\'e $x_i\in \overline \Z_p^*$.

\bigskip  On associe au morphisme
$\X$ un drapeau $(I_j)_{1\leq j \leq r}$ comme en (3.2).

 Le morphisme
$\tilde
\X$ tel que $\tilde \X (E_i)=x_i\phi (q)^{e_i}, \ x_i\in \overline
\Z_p^*, e_i \in \Q_{\geq 0},$ pour
$1\leq i
\leq n$ rel\`eve
$\X$ si et seulement si $(e_i,x_i)_{1\leq i \leq n}$ v\'erifient:
$$ e_{I_j}= t_j(t_j-1)/2 , \ \ \ r_p(x_{I_j})=\X(E_{I_j}) \leqno
 $$  pour tout
$1\leq j
\leq r$,
 et si
$I$ de cardinal $t$ n'appartient pas au drapeau
$(I_j)_{1\leq j \leq r}$, $$ e_{I }>t(t-1)/2, 
$$  
 o\`u l'on pose   $$x_I=\prod _{i\in I}x_i, \ \
e_I=\sum_{i\in I}e_i.$$ Le th\'eor\`eme 6 dit qu'il existe 
$(e_i,x_i)_{1\leq i
\leq n}$ satisfaisant ces conditions. Il suffit de  choisir  des
rel\`evements arbitraires $y_{I_{j}}\in
\overline \Z_p^*$ de
$\X(E_{I_j})$ pour 
$ 1\leq j \leq r$, et de prendre 

$e_i=(t_1-1)/2, \ x_i=y_{I_1}^{1/t_1}$, pour $i\in I_1$,  

$e_i=(t_{j+1}+t_j-1)/2, \ x_i=(y_{I_{j+1}}(y_{I_{j}}\ldots
y_{I_{j}})^{-1}) ^{1/(t_{j+1}-t_j)}$, pour
$i\in  I_{j+1}-I_j$ si $ 1\leq j \leq r-1$.

\noindent On v\'erifie imm\'ediatement les conditions sur les
$x_i$. V\'erifions les  conditions  sur les $e_I$. Soit $I\subset  
\{1, \ldots, n\}$ non vide avec $t=|I|$. Posons
 $x_j:=|I\cap I_j|$  pour
$1\leq j \leq r$. Donc $ t=\sum x_i$ et 
 $$2e_I=x_1 (t_1-1)  + (x_2-x_1)(t_2+t_1-1) +\ldots
+(x_r-x_{r-1})(n+x_{r-1}-1).$$
La suite $(x_j)$ est croissante et
  $x_j \leq t_j$ 
pour tout $j$, donc  
$$2e_I\leq x_1 (x_1-1)  +
(x_2-x_1)(x_2+x_1-1) +\ldots +(x_r-x_{r-1})(x_r+x_{r-1}-1).$$
avec  \'egalit\'e si et seulement si $I$ appartient au drapeau
$(I_j)_{1\leq j \leq r}$. Les propri\'et\'es voulues des $e_I$
r\'esultent
 de 
l'identit\'e  polyn\^omiale 
$$X(X-1)+ (Y-X)(Y+X-1)=Y(Y-1).$$
 qui implique que le membre de droite de l'in\'egalit\'e est
 $t(t-1) $.   

\bigskip {\bf 3.5 \ Preuve de la proposition  7 \ }  

Nous allons  d\'emontrer que la longueur d'un module standard supersingulier sur $\overline \F_{p}$ est $\geq 2^{n-2}$. Il suffit de  montrer qu'un caract\`ere supersingulier quelconque  de $A\cap H$ \`a valeurs dans $\overline \F_{p}$ poss\`ede   un rel\`evement $\X: A\cap H \to  \overline \Z_{p}$, tel que le module standard $I(\X, \overline \Q_{p})$ admette  $ 2^{n-2}$ sous-quotients simples et d'appliquer le th\'eor\`eme 5. La th\'eorie de Rogawski-Zelevinski permet de d\'ecrire les sous-quotients simples d'un module standard $I(\X, \overline \Q_{p})$ sur $ \overline \Q_{p}$ en fonction de $\X$ lorsque $\X(q)=q$. Maintenant que le principe est \'etabli, il reste \`a exhiber le rel\`evement.
Le rel\`evement d\'ecrit en 3.4 ne convient pas, mais sa forme est analogue \`a celui que l'on va exhiber. 

Soit $z\in \overline \F_{p}^{*}$ quelconque. On choisit deux unit\'es $a,b \in \overline \Z_{p}^{*}$ telles que $r_{p}(a^{n-1}b)=z$ et un nombre rationnel $u\in \Q$ tel que $0<u<1/n$. Un caract\`ere
$\X: A\cap H \to  \overline \Q_{p}$ tel que $\X(q)=q$ est d\'etermin\'e par les valeurs $\X(E_{i})$ pour tous les entiers $1\leq i \leq n$. On peut prendre ses valeurs quelconques. On prend
$$\X(E_{i})=aq^{u+i-1} \ \ \ (1\leq i \leq n-1), \ \ \ \X(E_{n})=bq^{(n-1)(1-u)}.$$
On v\'erifie d'abord que $\X$ est entier. Cela est \'equivalent \`a 
$$\X(E_{I})= q^{t(1-t)/2}\prod_{i\in I} \X(E_{i}) \in \overline \Z_{p},
 \ \ \ \X(E_{1,\ldots,n})\in  \overline \Z_{p}^{*},$$ 
pour toute partie $I$ de $\{1, \ldots, n\}$ d'apr\`es la description de $A\cap H$ par g\'en\'erateurs et relations (proposition 3.1). Ceci est  presque \'evident: on minore la somme $\sum_{i\in I}(u+i-1)$ par
$tu+ t(t-1)/2$ donc  
$$\X(E_{I}) \in a^{t}q^{tu+\N}$$
 est entier. On a  $\X(E_{1,\ldots,n})= a^{n-1}b$.  Donc le caract\`ere $\X$ est entier. Comme $u>0$ et $r_{p}(a^{n-1}b)=z$, la r\'eduction de $\X$ est le caract\`ere supersingulier \'egal \`a $z$ sur l'\'el'ement central $E_{1,\ldots,n}$.  Comme la longueur de $y_{i}$ ne d\'epend pas du choix de $i$, on voit que
  le caract\`ere de $X$ associ\'e est $\X$ contient un ``segment de Zelevinski'' de longueur $n-1$. 
  On d\'eduit de  que le module standard associ\'e contient $2^{n-2}$ sous-quotients irr\'eductibles distincts [Rog 5.2 page 456].

\bigskip {\bf Cas $n=2$ \ } La d\'ecomposition
des modules standards, donc la classification
des modules simples sur $\overline \F_p$  est connue [Vign\'eras].
On a  :

$\tilde \theta_1=q^{1/2}T_{w_o}T_{\w_1}^{-1},
\tilde \theta_2=q^{-1/2} T_{\w_1}$,  $E_i=q^{1/2}\tilde \theta_{y_i}$
pour
$ 1\leq i \leq 2$,  avec la relation
$E_1E_2=qZ$.

L'alg\`ebre commutative $A\cap H=\Z[q, \ Z^{\pm 1}, \  E_1,
E_2],$ 

Le centre   
 $Z(H)=\Z[q,\ Z^{\pm 1},Z_1=E_1+E_2]$. 

Les module standards sont de dimension $2$. 

Les modules standards supersinguliers sont
irr\'eductibles.

Les modules standards r\'eguliers  de
caract\`ere central $\w$ tel que $\w(Z_1)^2\neq \w(Z)$, sont
irr\'eductibles.

Les modules standards r\'eguliers  de
caract\`ere central $\w$ tel que $\w(Z_1)^2= \w(Z)$, sont
ind\'ecomposables de longueur $2$, contenant  caract\`ere trivial
(i.e. $T_s\to 0$)  et un caract\`ere signe
(i.e. $T_s\to -1$) comme sous-quotient.

Les modules simples se rel\`event.

\vfill\eject \centerline{ {\bf Bibliographie \ }}

\bigskip {\bf Bibliographie \ }

[Bourbaki AC] Bourbaki Nicolas \ {\sl Alg\`ebre commutative Chapitre 5 \`a
7}. Masson 1985.

[Bourbaki GAL] Bourbaki Nicolas \ {\sl  Groupes et alg\`ebres de Lie
Chapitres 4,5 et 6}. Hermann 1968.

[Cherednik] Cherednik Ivan \ {\sl A unification of
Knizhnik-Zamolodchikov and Dunkl operators via affine Hecke algebras}.
Invent. math. 106, 411-431 (1991).

[Haines] Haines Thomas J. \ {\sl The combinatorics of Bernstein
functions.} 
    Trans. Am. Math. Soc. 353, No.3 (2001), 1251-1278.

[Lusztig 1] Lusztig Georges \ {\sl Affine Hecke algebras and their
graded version.} Journal of A.M.S. Vol. 2, No 3, 1989.

[Lusztig 2] \ \ \  {\sl Singularities, character formulas, and a
$q$-analog of weight multiplicities.} Ast\'erisque 101-102 page
208-229. SMF 1984.

[Lusztig 3] \ \ \ {\sl Hecke algebras with unequal parameters.}
ArXiv:math.RT/0208154 2002.

[Ollivier] Ollivier Rachel \  {\sl Modules simples en caract\'eristique
$p$ des alg\`ebres de Hecke affines de type $A_1, A_2$. }DEA Institut
de Math\'ematiques de Jussieu, 28 Juin 2002.

[Rogawski] Rogawski J. D. \ {\sl On modules over the Hecke algebra of
a
$p$-adic group}. Invent. math. 79, 443-465 (1985).

[Steinberg] Steinberg R. \ {\sl On a theorem of Pittie}. Topology 14
(1975) 173-177.

[Vign\'eras] Vign\'eras Marie-France \ {\sl Representations modulo $p$
of the
$p$-adic group $GL(2,F)$}. Pr\'epublication 301, Institut de
Math\'ematiques de Jussieu, Septembre 2001. A paraitre \`a Compositio
Mathematicae.

\vfill\eject

\centerline {{\bf APPENDICE }}

\bigskip {\bf D\'efinitions} [Springer T.A.
Reductive groups. Proceedings of Symposia in Pure Mathematics. Vol. 33
(1979), part 1, pp. 3-27] \   Une donn\'ee radicielle bas\'ee
 est un
quintuplet
$(X,X^\vee,R,R^{\vee},B)$ o\`u  

-  $X,X^\vee$ sont des
groupes libres ab\'eliens de rang fini $\geq 1$, en forte dualit\'e
par  un forme bilin\'eaire  
$(\ , \ ):X\times X^\vee \to \Z$,
$R$ et
$R^\vee$ sont des sous-ensembles finis non vides de $X$ et de $X^\vee$ 
dont les
\'el\'ements sont appel\'es ``racines'' et   ``coracines'', $B$ est
une base de $R$.

- Une bijection 
$\alpha \leftrightarrow  \alpha^\vee$ est donn\'ee entre les racines
 et les coracines  telle que 
 $(\alpha, \alpha^\vee)=2$,

-  
les
reflexions $s_{\alpha}$ dans $GL(X)$ et dans $GL(X^\vee)$
  d\'efinies
par
$s_{\alpha}(x)=x-(x,\alpha^\vee)\alpha$ et 
$s_{\alpha}(x ^\vee)=x^\vee-( \alpha,x^\vee)\alpha^\vee$  
pour tout $\alpha \in  R, \alpha^\vee\in R^\vee, x\in X,
x^\vee
\in X^\vee$, permutent les racines $s_{\alpha}(R)=
R$,  et 
  les coracines $s_{\alpha^\vee}(R^\vee)= R^\vee$.

On supposera  le syst\`eme de racines $R$  r\'eduit.

\bigskip  Le groupe de Weyl fini $W_o$ est le sous-groupe de $GL(X)$  
engendr\'e par les r\'eflexions $s_{\alpha}, \ \alpha \in R$. Le
syst\`eme $(W_o, \{ s_{\alpha}, \ \alpha\in B \})$ est un syst\`eme de
Coxeter. Il est isomorphe canoniquement au sous-groupe de $GL(X^\vee)$
engendr\'e par les  r\'eflexions $s_{\alpha^\vee}, \ \alpha^\vee \in
R^\vee$, et la forme bilin\'eaire $(\ , \ )$ sur $X\times X^\vee$ est
$W_o$-invariante.

Le groupe de Weyl $W$ est le produit
semi-direct
  de $W_o$ et de $X$.  
On notera
$e^x$ l'\'el\'ement
$x\in X$ plong\'e dans $W$; on \'ecrira  $W=W_o.e^X$. On a
$w_oe^x=e^{w_o(x)}w_o$ pour $(w_o,x)\in W_o\times X$.

Le groupe de Weyl affine $W_{aff}$ est le produit
semidirect
 de $W_o$ et du sous-groupe $Q(R)$ de $X$
engendr\'e par
$R$; on  \'ecrira $W_{aff}=W_o e^{Q(R)}$.

On d\'efinit un ordre partiel $\preceq$ dans $R^\vee$ tel que
$\alpha_1^\vee\preceq \alpha_2^\vee$ si $\alpha_2^\vee -\alpha_1^\vee=
\sum _{\beta
\in B}n_{\beta}\beta^\vee$ pour des entiers $n_{\beta}\geq 0$. On note
$R_m$ les racines  $\alpha \in R$ telles que les coracines
$\alpha^\vee \in R^\vee$ sont minimales  
 pour l'ordre
$\preceq$ et 
$$S=\{ s_\beta, \ \beta\in B\} \cup \{ e^{-\alpha}s_{\alpha}, \
\alpha \in  R_m\}.$$ 
  Le syst\`eme $(W_{aff},S)$ est un syst\`eme de Coxeter. 
 
 Toute racine est combinaison
lin\'eaire
\`a coefficients entiers de m\^eme signe d'\'el\'ements de $B$. Si le
signe est $\geq 0$  la racine est dite positive,
 et n\'egative sinon. On note $R^+$  l'ensemble des racines
positives et $R^-$ celui des racines n\'egatives. 
Un
\'el\'ement
$x\in X$ tel que
$(x,\beta^{\vee})\geq 0$ pour tout $\beta\in  B$ est dit
dominant. 

 Le monoide $X_{dom}$ des \'el\'ements dominants de $X$ est 
de type fini.

\bigskip  {\bf   Longueur \ } 
  La  longueur 
$\ell$    de $W=\Omega.W_{aff} $   prolonge celle du syst\`eme de
Coxeter $(W_{aff}, S)$ [Bourbaki GAL
IV.1.1 page 9], de sorte que la longueur est
constante sur $\Omega w_{aff}\Omega$ pour tout
$ w_{aff}\in W_{aff}$, ce qui a un sens car $\Omega$
normalise   $S$. La longueur   v\'erifie les
propri\'et\'es:

- Le groupe
$\Omega$ est form\'e par les \'el\'ements   de longueur
$0$.

- La longueur est invariante par passage
\`a l'inverse $\ell(w)=\ell(w^{-1})$ pour tout $w\in W$.

- La longueur de $W=W_oe^X$ s'exprime par une somme d'entiers
naturels  index\'ee par les racines positives  
  [Iwahori-Matsumoto IHES 25, section I.10], g\'en\'eralisant celle  
pour le groupe de Weyl fini $W_o$  [Bourbaki GAL VI 1.6 page 157].  Pour
$w_o\in W_o, x\in X$, 

$$\ell(w_o e^x)=\sum_{\alpha \in R ^+, w_o(\alpha)\in
R^+}|(x,\alpha^{\vee})| +\sum_{\alpha \in R ^+, w_o(\alpha)\in
R^-}|1+(x,\alpha^{\vee})|.$$

- La longueur de $w_o\in W_o$ est donc le nombre de racines 
positives $\alpha \in R^+$ telles que $w_o(\alpha)\in R^-$. 

- La longueur sur $e^X$ est $W_o$-invariante:
$$\ell(w_o(e^x))=\ell(e^x) ,  \ \ \ \ w_o\in W_o, x\in X.$$
 En effet, il suffit de le
v\'erifier pour une r\'eflexion $s_{\beta}$ d\'efinie par une racine simple  $\beta
\in B$, et dans ce cas d'utiliser  que:

a)  $(s_{\beta}(x),\alpha^\vee)=(x ,s_{\beta }(\alpha)^\vee ) $ pour
$\alpha\in R $,

b) la r\'eflexion
$s_{\beta}$ permute les  racines positives
 diff\'erentes de 
$\beta$  [Bourbaki GAL VI   1.6 Cor 1 page 157].

- La longueur de $e^x$ pour $x\in X_{dom}$  est
un entier pair, car $\ell(e^x)=(x,2\rho^\vee)$ o\`u
$$ 2\rho^\vee:=\sum_{\alpha \in R^+}
\alpha^\vee  = 2\sum _{\beta \in B} \w_{\beta^\vee}$$ 
est deux fois la somme des  poids fondamentaux $ \w_{\beta^\vee} $
tels que $$(\alpha, \w_{\beta^\vee}) =\delta_{\alpha,\beta} 
, \ \ \ \alpha, \beta \in B$$ [Bourbaki GAL VI 1.10  prop. 29 page
168]. En particulier,

- $\ell(w_oe^x)=\ell(w_o)+\ell(e^x)$  si $x$ est dominant.  

 - Soient $x,x'\in X$ et $w_o\in W_o$. Posons $n(\alpha,w_oe^x)=
 (\alpha^\vee, x)$ si $w_o(\alpha)\in R^+$   et $n(\alpha,w_oe^x)= 
 1+(x,\alpha^\vee )$ si
$w_o(\alpha)\in R^-$. 
 Si les entiers
$n(\alpha,w_oe^x)$ et
$n(\alpha,e^{x'})$ ont le m\^eme signe au sens large pour tout $\alpha
\in R^+$, alors
$$\ell (w_oe^{x+x'})=\ell(w_oe^x)+\ell
(e^{x'}).$$ 

En particulier,
   
- $\ell (w_oe^{2x})=\ell(w_oe^x)+\ell
(e^{x})$.

 - $\ell(e^{x+x'})=\ell(e^x)+\ell(e^{x'})$ pour tout
$w_o\in W_o$ et tout $x,x'\in w_o(X_{dom})$.

\bigskip {\bf (Ap.1) \ } On v\'erifie:

 {\sl Soient $u,v\in W_o$ et $x\in X$. Alors,

 a) $\ell(u)+\ell(v)-\ell(uv) $ est deux fois le nombre de racines $\alpha\in
R^+$ telles que

\centerline { $ v(\alpha)\in R^-, \  uv(\alpha)\in R^+$.}

b)  $\ell(u)+\ell(ve^x)- \ell(uve^x)$ est deux fois le nombre  de racines $\alpha\in R^+$  
telles que  
  
\centerline {1) $ v(\alpha)\in R^-, \ uv(\alpha)\in R^+, \ (x,\alpha^\vee)\geq 0$,}

ou 

\centerline { 2)
$ v(\alpha)\in R^+, \ uv(\alpha)\in R^-, \ (x,\alpha^\vee)<0$.}}

\bigskip En particulier, 

- si $\alpha, \beta\in B$, alors 
 $\ell
(s_{\alpha}e^{\beta})=1$ si et seulement si $(\beta, \alpha^\vee)<0$.

\bigskip {\bf (Ap.2) \ }  On va d\'emontrer:

{\sl 

Soit $x\in X$. Notons  $\ell(x)$ le nombre  de  racines
positives $\alpha\in R^+$ telles que
 $(x,\alpha^\vee)<0$. 

1) La longueur d'un \'el\'ement $w_o
\in W_o$  tel que $ w_o(x) $ est dominant est $\geq \ell(x)$. 

Il existe un unique $u 
\in W_o$ de longueur  $\ell(x)$  tel que $y= u (x) $ est dominant.

On a
$(x,\alpha^\vee)<0$ si et seulement si 
 $u_o(\alpha)\in R^-$, pour tout   $\alpha\in R^+$.

2) Pour tout $t$ tel que $1\leq t
\leq \ell(x)$, soit  $\beta_t\in 
B$ une racine  simple  de r\'eflexion associ\'ee $s_t$, telle  que $u=
s_{
\ell(x) }
\ldots s_2s_1$ soit une d\'ecomposition r\'eduite; 
 posons  $u_t:=s_t\ldots
s_2s_1 $ et $u_o:=1$.

Pour tout $w_o\in W_o$, on a
$$ \ell(w_o )=    \ell(w_ou_t^{-1})+\ell ( u_{t } ) 
 $$
si $w_o(\alpha)\in R^-$ pour tout $\alpha \in \{ \beta_1,\ 
s_1(\beta_2),\  \ldots ,  \ s_1\ldots
s_{ t-1 }(\beta_{t}) \} $, et
$$ 
\ell(w_oe^x)=    \ell(w_ou_t^{-1})+\ell ( u_{t } e^x) 
$$
si $w_o(\alpha)\in R^+$ pour  tout
$\alpha \in \{s_1\ldots
s_t(\beta_{t+1}) ,\     \ldots ,  \ s_1\ldots
s_{ \ell(x)-1} (\beta_{\ell(x)}) \}. $  

En particulier,
  $\ell(ue^x)=\ell(e^x)-\ell (x)$.}

\bigskip Preuve.  1)  Soit $(w_o,x)\in W_o\times X$. On a
$(w_o(x), w_o(\alpha)^\vee)=( x , \alpha ^\vee)$ pour toute racine
  $\alpha \in R $.   On voit que {\sl $w_o(x)$
est dominant, i.e. $\ell (w_o(x))=0$,  si et et seulement si
$w_o (\alpha)\in R^-$ est une racine n\'egative pour  toute racine
positive
$\alpha
\in R^+$ telle  que $(x,\alpha^\vee)<0$.}
Donc si $w_o(x)$ est dominant, la longueur de $w_o$ est   $\geq \ell
(x)$ (on rappelle que $\ell(x)$ n'est pas la longueur de $e^x$).
 
Supposons que $x $ n'est pas dominant, i.e. $\ell(x)> 0$. On construit un
\'el\'ement $u\in W_o$ de longueur minimale $\ell(x)$ tel que $u(x)$ est dominant de la fa\c{c}on suivante.  

Il existe au moins une  racine
simple
$\beta $ telle que
$(x,\beta ^\vee)<0$.   

On a  
$\ell(s_{\beta}(x))=\ell(x) -1 $ car
 $(s_{\beta}(x),\alpha^\vee)=(x, s_{\beta}(\alpha^\vee))$ et $s_{\beta}$ permute les racines
positives diff\'erentes de
$\beta$.

 Posons $\beta_1=\beta, u_1=s_{\beta_1}, x_1=s_{\beta}(x) $. 

On recommence    en
partant de
$x_1$. Au bout de $\ell(x)$
 \'etapes on obtient un \'el\'ement dominant. On a   ainsi choisi une suite de
$\ell(x)$
   racines simples     
$\beta_t $, d'\'el\'ements $u_t=s_{\beta_t }\ldots
s_{\beta_2}s_{\beta_1}$ de $W_o$, d'\'el\'ements $x_t=u_t(x) $ de
$X$, tels que 
$$  
  (x_{t-1} , \beta_t^\vee)<0, \ \ \ \ell(x_t)=\ell(x)-t$$  pour
tout
$1\leq t \leq \ell(x)$. On prend $u=u_{\ell(x)}$. 
 
Les racines $\alpha \in R^+$ telles que $u(\alpha)\in R^-$ sont  
[Bourbaki GAL VI,
\S 1, 1.6 Cor. 2] :
$$  \beta_1,\  s_{\beta_1}(\beta_2),\  \ldots ,  \ s_{\beta_1 }\ldots
s_{\beta_{\ell(x)-1}}(\beta_{\ell(x)}). \leqno (S) $$
Ce sont les $\ell(x)$ racines telles que 
$(x,\alpha^\vee)<0$ car $(x_{t-1} , \beta_t^\vee)=(x,
s_{\beta_1} \ldots s_{\beta_{t-1}}(\beta_t)^\vee) $.

L'unicit\'e de $u$ provient de ce que l'ensemble des racines positives $\alpha \in R^+$ telles que $u(\alpha)\in
R^-$ est ind\'ependant de $u$. Cet ensemble d\'etermine $u$
(introduire l'ensemble not\'e $T_u$ comme dans [Bourbaki GAL VI
 1.4 page 13-14] qui d\'etermine $u$, puis appliquer [Bourbaki GAL VI 1.6
prop.17 page 157]).

2)  Soit $w_o\in W_o$. Si $t\geq 1$, les racines $\alpha \in R^+$ telles que $u_t(\alpha)\in R^-$ sont les $t$
premiers termes de la suite (S). On applique (Ap.1.1) pour obtenir les
deux premi\`eres \'egalit\'es sur la longueur.

Le cas particulier s'obtient en   prenant $w_o=1, t=\ell(x)$ dans la
seconde \'egalit\'e.

\bigskip {\bf Ordre de Chevalley-Bruhat sur $W$ \ } C'est un ordre partiel $\leq $ sur $W$ qui prolonge l'ordre de
Chevalley-Bruhat usuel sur le syst\`eme de Coxeter $(W_{aff}, S)$.  
Les propri\'et\'es suivantes pour
$w,w'\in W_{aff}$ sont \'equivalentes [Lusztig3], voir aussi
le lemme 8.11 dans [Bernstein I.N., Gelfand I.M. and
Gelfand S.I., Schubert cells and
cohomology of the spaces $G/P$,  Russ. Math. Surv. 28 (1973), 1-26] :

- il existe une expression r\'eduite de $w$ telle qu'en omettant certains termes on obtient une expression
  de $w'$,

- pour toute  expression r\'eduite de $w$,   en omettant certains termes on obtient une expression
  de $w'$,

- il existe une suite d'\'el\'ements $w_o=w', w_1, \ldots, w_k=w$ dans $W_{aff}$ telle que 

\centerline {$\ell(w_1)-\ell(w_o)=\ldots = \ell(w_k)- \ell(w_{k-1}^{-1})=1, \ \ \ w_1w_o^{-1}, \ldots,
w_kw_{k-1}^{-1}\in T,$} 

- il existe une suite d'\'el\'ements $w_o=w', w_1, \ldots, w_k=w$ dans $W_{aff}$ telle que 

\centerline {$\ell(w_1)-\ell(w_o)=\ldots = \ell(w_k)>\ell(w_{k-1}) , \
\ \ w_1w_o^{-1}, \ldots, w_kw_{k-1}^{-1}\in T$.} 

\noindent o\`u $T$ est l'ensemble  des conjugu\'es de $S$ dans
 $W_{aff}$.

Par d\'efinition de l'ordre de Chevalley-Bruhat, $w'\leq w $ si  ces conditions sont r\'ealis\'ees.

 La propri\'et\'e
d'\'echange des syst\`emes de Coxeter [Bourbaki GAL IV, \S 1, 1.5 page 15]
montre que dans $W_{aff}$:

\centerline {  $\ell(w )+\ell(w' )=\ell(w w' )$ implique
     $w \leq w w' $. }

  L'ordre de Chevalley-Bruhat de $W_{aff}$ est
invariant par conjugaison par $\Omega$, car $\Omega$ normalise
 $S$.

On prolonge naturellement l'ordre de Chevalley-Bruhat de $W_{aff}$ \`a
$W=\Omega  W_{aff}$. Par d\'efinition $u w_{aff}  \leq u'w'_{aff} $ si et seulement si $u=u', 
w_{aff} \leq w'_{aff} $, pour tout $ u,u' \in \Omega, w_{aff}, w'_{aff} \in W_{aff}$.
On a les propri\'et\'es
  \'equivalentes suivantes:

- $w_{aff} \leq w'_{aff} $ 
 si et seulement si $w_{aff}u\leq w'_{aff}u$ pour tout $u\in \Omega,
    w_{aff}, w'_{aff} \in W_{aff}$,

- $w_{aff} \leq w'_{aff} $ 
 si et seulement si $uw_{aff}  \leq uw'_{aff} $ pour tout $u\in \Omega, w_{aff}, w'_{aff} \in W_{aff}$, 

- $w_{aff} \leq w'_{aff} $ 
 si et seulement si $u
 w'_{aff}u'\leq u w_{aff}u'$  pour tout $ u,u' \in \Omega, w_{aff}, w'_{aff} \in W_{aff}$.

On v\'erifie que:

Si $w,w'\
\in W, w'\leq w $, alors ${w'}^{-1}\leq w ^{-1}$ et $\ell(w')\leq \ell(w)$.

\end